\pretolerance=2000
\hbadness=2000
\vbadness=2000

\parindent=1truecm

\def\id{1\kern-2.5pt{\rm l}}
\def\dd{\partial}

\def\tr{{\rm Tr}}

\def\so{\hbox{\got so}}

\font\got=eufm10

\font\ccap=cmbx12 at 18pt
\font\scap=cmbx12 at 12pt

\font\sscap=cmbxsl10

\def\S{\hbox{\eusm S}}

\def\F{\hbox{\got F}}
\def\M{\hbox{\got M}}

\def\ad{\rm ad}
\def\l{\gamma}
\def\ga#1#2{\gamma_{_{#1 #2}}}

\def\S{\hbox{\got S}}
\def\H{\hbox{\got H}}

\def\pd#1#2{#2_{_{#1}}}

\def\var#1#2{{\delta#1\over{\delta#2}}}

\def\pois#1{\lbrace#1\rbrace}

\vskip 7truecm

\centerline{\ccap On the Hamiltonian and Lagrangian}
\medskip
\centerline{\ccap structures of time--dependent reductions }
\medskip
\centerline{\ccap  of evolutionary PDEs}
\bigskip
\vskip 3truecm
\centerline{\bf Monica Ugaglia}
\centerline{\bf Via Aldo Moro 37, 57024 Donoratico, Italy}
\centerline{\bf e-mail: ugaglia@grog.scumm.it}
\vskip 5truecm
\noindent{\sscap Abstract:} In this paper we study the reductions of 
evolutionary PDEs on the manifold of the stationary
points of  time--dependent symmetries. In particular we describe how that the 
finite dimensional 
Hamiltonian structure of the reduced system is obtained from the Hamiltonian 
structure of the initial PDE and we construct the time--dependent Hamiltonian
function. We also present a very general Lagrangian formulation of 
the procedure  of reduction. As an application we consider the case of the 
Painlev\'e equations PI, PII, PIII, PVI and also certain higher order systems
appeared in the theory of Frobenius manifolds.
\vskip 5truecm
\centerline{Preprint SISSA 11/99/FM}
\centerline{To appear in {\it Differential Geometry and its
Applications}}
\vfill\eject

\noindent{\scap  Introduction}
\bigskip

Consider an evolutionary 
PDE with one spatial variable
$$
u_t=F(u,u_x,\ldots,u^{(m)})\eqno{(0.1)}
$$
 and a symmetry of this equation,
 i.e. another system
$$
u_s=G(u,u_x,\ldots,u^{(k)})\eqno{(0.2)}
$$
commuting with the first one,
$$
(u_s)_t=(u_t)_s.
$$
 The set of the stationary points $u_s=0$  of the symmetry is a 
finite-dimensional invariant
 manifold for the system (0.1). Particularly, in important examples, the invariant manifold can be described as a set of stationary points of a first integral of the system (0.1):
$$ 
\var I {u(x)}=0,
$$

$$ 
I=\int L(u,u_x,\ldots,u^{(n)})dx,\qquad {dI\over dt}\equiv 0.
$$
In this case it is known ([BN], [Mo]) that 
the restriction of the initial PDE  to the
 invariant submanifold  is a Hamiltonian system of ODEs.
In particular Bogoyavlenskii and Novikov [BN] found a universal scheme 
to construct the Hamiltonian function 
 of the
reduced system in terms of the Hamiltonian of the original PDE.
In this paper we  extend this scheme to more general finite dimensional
invariant
 submanifolds specified by local $x$- and time-dependent symmetries 
and conservative quantities of the evolutionary equation.
To distinguish this class of symmetries from the previous one we will 
call them {\it scaling symmetries}. We show that the
restriction of the starting equation  on the finite dimensional manifold
admits a natural description as a Hamiltonian system 
with time--dependent Hamiltonian.

The best known class of examples of evolutionary PDEs admitting 
nontrivial symmetries and conservation laws are integrable systems of 
soliton theory (see [SM] and references therein).
The finite dimensional manifolds of the stationary points of integrable systems are typically described
by ODEs of Painlev\'e type [AS],[CD]. For the simplest
examples of these restrictions the Hamiltonian structure is already
known. For example for the classical six Painlev\'e equations 
the Hamiltonian description was found by Okamoto, [O].
Although the relationship between the starting 
PDE and the reduced ODE is clear and has been investigated 
quite a lot (see, e.g. [AS]), the relationship
 between the starting Hamiltonian structure and the 
reduced one has not been elucidated. 
This work will give a contribution in understanding of this relationship. 

As a  first result (see section 2) 
we prove that the finite dimensional Hamiltonian
structure of the  ODEs is obtained from the Hamiltonian structure of the
starting 
PDE, via scaling reduction. Particularly, we
 construct the
time--dependent Hamiltonian function of the reduced system. In the
time--independent case this procedure coincides with the
 well known stationary--flow reduction discovered by Bogoyavlenskii and Novikov
[BN]. As an application we present the case of PI, PII, PIII, PVI and also
certain higher order systems appeared recently in the theory of Frobenius
manifolds [D1].   

As a second result (see section 3) we present a very general
Lagrangian formulation of the procedure of reduction of an evolutionary 
system (0.1).
Namely, we prove that this restriction is again a Lagrangian system with the
Lagrangian function $\Lambda$, such that
$$
{d\Lambda\over dx}={d L\over dt}.
$$
\bigskip

\noindent The work is structured as follows: after recalling, in Section 1, 
some basic 
facts about the Hamiltonian structure
of the evolutionary PDEs, and  briefly summarizing 
the   method of reduction of evolutionary flows on the
manifold of stationary points of their integral, introduced by Bogoyavlenskii
and Novikov [BN], in Section 2 we consider the generalization
of this procedure to  scaling symmetries. The reduced flow is a 
time--dependent Hamiltonian system, and in 
Theorem 2.1 we give the relationship between  the infinite--dimensional
Hamiltonian structure and the reduced one.

Section 3 is devoted to a Lagrangian approach to the problem: after 
describing the general framework, in Theorem 3.1 we give the procedure
of reduction and we construct the reduced Lagrangian function. In Section 3.2 
we establish
the relationship with the Hamiltonian approach. As an application we study 
the Lagrangian reduction of KdV on the manifold of the fixed points of the 
$7$--th flow.

Section 4 contains the application of the theory to the scaling
reductions from KdV,  mKdV and Sine--Gordon equations respectively  
to Painlev\' e I, Painlev\' e II and III. These examples are studied both from 
the Hamiltonian
 and  the Lagrangian point of view. 

In Section 5 we study the $n$-waves equation and his scaling 
reduction to a system 
of commuting Hamiltonian flows on the Lie algebra $\so(n)$. 
 The reduced system
is a  non-autonomous Hamiltonian
system w.r.t. the Poisson structure of $\so(n)$.
In particular, for $n=3$, adding an additional symmetry condition,
one arrives at  Painlev\' e VI equation.
\bigskip
\bigskip
\noindent {\scap 1. Infinite dimensional Hamiltonian structures and 
stationary flows reduction}
\bigskip
\noindent  Let us consider the phase space $\M$ of smooth maps of
the circle into some smooth $n$-dimensional manifold. Actually we can forget
about the boundary conditions when dealing  with local functionals only.
We denote by $\F$ the space of  smooth functionals on $\M$ of the form $$
F(u)=\int f(x,u(x),u_x(x),\ldots ,u^{(m)}(x))dx, $$ where the density $f$ depends
only on a finite number of derivatives of $u$. On the space $\F$ the variational
derivative  $\var F {u^i(x)}$ is defined by $$ \delta F=\int \var F
{u^i(x)}\delta u^i(x) dx. $$ Explicitly,
$$
\var F {u^i(x)}={\partial f\over
{\partial {u^i(x)}}} +\sum (-1)^{^k}
{d^{^k}\over {dx^{^k}}} {\partial f\over {\partial {u^{{^i}(k)}(x)}}}.
$$
One can define on $\M$ the (formal) Poisson brackets
$$
\pois{u^i(x),u^j(y)}=w^{^{ij}}(x,y)=\sum^N_{k=0}A^{^{ij}}_k
\delta^{^{(k)}}(x-y),
$$
 where $A^{^{ij}}_k$ depends on a finite
number of derivatives of $u$. 
This induces on $\F$ the Poisson bracket
$$
\pois{F,G}=\int \var F {u^i(x)} P^{^{ij}} \var G
{u^j(x)}dx,
$$
 where
$$
P^{^{ij}}=\sum^N_{k=0}A^{^{ij}}_k({d\over dx})^{^{k}}.
$$

\noindent A Hamiltonian system on $\M$ has then the form 
$$ 
u^i_t(x)=\pois{u^i(x),H}= P^{^{ij}} \var H{u^j(x)}. 
$$
\bigskip
\noindent In particular, we consider so called Gardner--Zakharov--Faddev bracket
$P^{^{ij}}=\delta ^{ij}{d\over dx}$. In this case a Hamiltonian system has 
the form
 $$
u_t(x)=\pois{u(x),H}={d\over dx} \var H{u(x)}\eqno{(1.1)}
$$
with Poisson bracket
$$
\pois{F,G}=\int \var F {u(x)}{d\over dx} \var G
{u(x)}dx.
$$
Let us consider a first integral 
$$
I=\int L(x,u(x),u_x(x),\ldots ,u^{(n)}) dx,
$$
where $L$ does not depend on $t$. 
The
 generalized Euler--Lagrange
 equation
$$
\var I {u(x)}=0\eqno{(1.2)}
$$
generically
is a ODE of order $2n$ fixing the $2n$--dimensional manifold  $\S$ of the
stationary points of the first integral $I$. Because of the Lax lemma (see 
[Mo]) this submanifold 
is invariant w.r.t the evolutionary equation (1.1). The functional $L$ is the
Lagrangian of the $x$--flow defined by (1.2). If $L$ is nondegenerate, then it
defines also on $\S$ the natural system of canonical coordinates 
$$
\eqalign{q_{_i}&=u^{(i-1)},\ \ \ \ i=1,2,....,n\cr 
p_{_i}&=\var  I  {u^{(i)}}.\cr}
$$ 
and  equation (1.2) can be put in the Hamiltonian form
$$
\cases{
(p_{_i})_x=-{\dd { H}\over{\dd q_{_i}}}\cr
(q_{_i})_x={\dd { H}\over{\dd p_{_i}}},\cr}
$$
where $H$ is the generalized Legendre transform of $L$:
$$
H=-L+\sum ^{n}_{i=1} \var { I}  {u^{(i)}}
u^{(i)}
$$
which, in terms of the canonical coordinates takes the form:
$$
 H= - { L}
+\sum ^{n}_{1}p_{_i}{dq_{_{i}}\over dx}.
$$
It is well known that the starting PDE can be restricted on $\S$ and the 
restriction is a Hamiltonian system of ODEs. In particular
 Bogoyavlenskii and Novikov discovered the algorithm to construct
the Hamiltonian functions of the
reductions in terms of the Hamiltonian of the original evolutionary equation.
They considered the case of a hierarchy of evolutionary equations
$$
{d u\over dt_{_{k}}}={d\over dx} \var {\pd k I}  {u(x)}
$$
with $\pd k I=\int \pd k L (u(x),u_x(x),\ldots ,u^{(n_k)}(x))dx$,
and they described the  reduction procedure of the $k$--th flow   on the
finite dimensional manifold of the stationary points of the $j$--th flow.
They  proved that all the flows of the hierarchy
reduce to 
finite dimensional Hamiltonian system. The  Hamiltonian
function for the reduced $k$--th flow, $(-Q_{k,j})$, is determined by :

$$
 \var {\pd j I} {u(x)} {d\over dx} \var {\pd k I} {u(x)} \equiv
 {d\over dx} Q_{k,j}.
$$
Mokhov [Mo] generalized this result to not necessarily Hamiltonian evolutionary
PDEs.
\bigskip

\bigskip
\noindent {\scap 2. Scaling reductions of evolutionary 
systems:}

\noindent {\scap Hamiltonian formulation}
\bigskip
\bigskip

\noindent In this Section we  extend the
Bogojavlenskii--Novikov scheme to  finite dimensional invariant
 submanifolds specified by time--dependent local
symmetries.

We start from  a partial differential equation of order $m$ on the
functional space $\M$, describing the evolution of the function $u(x)$ in the
time $t$ and a scaling symmetry
$$
u_s=G(x,t,u,u_x,\ldots,u^{(k)}).
$$
Our main assumption is that the set of stationary points of the symmetry
can be formally represented in the Euler--Lagrange form
$$
\var I {u(x)}=0,
$$

$$ 
I=\int L(x,t,u,u_x,\ldots,u^{(n)})dx,\qquad {dI\over dt}\equiv 0.
$$ 
 It is an ordinary differential equation of order
$2n$ depending explicitly on the parameter $t$. If $L$ is
nondegenerate, the space of the solutions is a $2n$ dimensional manifold  $\S$,
which naturally carries a  system of canonical coordinates. As in Section 1 we 
will show that, in these
coordinates,  the Euler--Lagrange equation is Hamiltonian, with Hamiltonian function $H$,
obtained from $L$ via Legendre transform:
$$
H=-L +\sum ^{n}_{i=1}p_i {d q_i\over dx}.
 $$
 Following the scheme of [BN], we prove that one
can reduce on $\S$ also the equation of the evolution in $t$, which 
results to be a
Hamiltonian system. We also give a universal scheme to produce the time-dependent
 Hamiltonian function of this reduced system. Indeed the following theorem holds:
\bigskip

\noindent {\bf Theorem 2.1:}  {\it If the evolutionary PDE:
$$
u_t=F(u,u_x,\ldots,u^{(m)}),
$$
admits a   nondegenerate {\it scaling} symmetry,
then,
 on the manifold $\S$ of the stationary points of the symmetry:
$$
\var I {u(x)}=0,
$$

$$
I=\int L(x,t,u,u_x,\ldots,u^{^{(n)}})
 dx,\qquad {dI\over dt}\equiv 0,
$$
 it reduces to a Hamiltonian motion in $t$, for the
time dependent Hamiltonian function $(-\tilde Q)$, that is the
reduction on $\S$ of $$
Q=\Lambda -\sum ^{n}_{i=1}p_i {d q_i\over dt},
\eqno{(2.1)}
$$
 where $p_i,q_i$ are the canonical coordinates on $\S$, expressed in terms
of $u,u_x,\ldots,u^{^{(2n-1)}}$,  and the function $\Lambda$ is determined by}
  $$
{dL\over dt}={d\Lambda\over dx}.\eqno{(2.2)}
$$
 
\bigskip
\noindent {\bf Proof:} We prove the theorem in three steps: first we describe
the submanifold $\S$ of stationary points of the symmetry $I$, where we
introduce a system of canonical coordinates; then we deduce, on $\S$,  a {\it
zero--curvature} equation for $(- \tilde Q)$ and the Hamiltonian function
$H$ of the reduced $x$--flow. Finally
 we prove that  the restricted $t$--flow is Hamiltonian
on $\S$,with Hamiltonian function $(-\tilde Q)$.
 \bigskip
1)  The manifold $\S$  is the $2n$--dimensional manifold of the
solutions of the Euler--Lagrange equation
$$
\var  {I}  {u(x)}=0.\eqno{(2.3)}
$$

\noindent It is invariant under the $t$--flow and it naturally carries a
system ofcanonical coordinates:
 $$
\eqalignno{q_{_i}&=u^{(i-1)},\ \ \ \ i=1,2,\ldots,n&(2.4 a)\cr 
p_{_i}&=\var  I  {u^{(i)}},&(2.4 b)\cr}
$$ 
obtained via generalized Lagrange transform (here we suppose that the
generalized Lagrangian $L$ is nondegenerate). Observe that now the $p_i$
 depend on $x$ and on  $t$.

Reversing relations (2.4), one can express the derivatives  
$u,u_x,\ldots,u^{(2n-1)}$ in terms of the
canonical coordinates $p_i$ and $q_i$, $x$ and $t$; explicitely:
$$
\cases{u^{(n)}=(q_{_{n}})_x=g_1(x,t,q_{_1},\dots,q_{n},p_{n})\cr
u^{(n+1)}=g_2(x,t,q_{_1},\dots,q_{n},p_{n},p_{n-1})\cr
\dots \dots\cr
u^{(2n-1)}=g_n(x,t,q_{_1},\dots,q_{n},p_{n},\dots,p_{_1}).\cr}
$$
Observe that (2.4) gives the identities:
$$
\cases{(p_{_1})_x+{\dd  H\over{\dd q_{_1}}}\equiv -
\var { I}  {u}\cr
(p_{_i})_x+{\dd { H}\over{\dd q_{_i}}}\equiv 0,\qquad i>1\cr
(q_{_i})_x-{\dd { H}\over{\dd p_{_i}}}\equiv 0,\cr}\eqno(2.5)
$$
where $H$ is the generalized Legendre transform of $L$:
$$
-L+\sum ^{n}_{i=1} \var  I  {u^{(i)}}
u^{(i)}
$$
which, in terms of the canonical coordinates takes the form:
$$
 H= -  L
+\sum ^{n}_{1}p_{_i}{dq_{_{i}}\over dx}.
$$

The first of  identities (2.5) allows us to express the higher derivatives
$u^{(m)}$ for $m\geq 2n$ in terms of $x$,$t$, $p_i,\ q_i$ and
$p_{{_{1}}}^{{_{(l)}}}$ with $l=1,\dots,{m-2n+1}$, explicitly:
$$
\cases{u^{(2n)}=g_{n+1}(x,t,q_{_1},\dots,q_{n},p_{n},\dots,p_{_1},
(p_{_1})_x)\cr
\dots\dots\cr
u^{^{(m)}}=g_{m-n+1}(x,t,q_{_1},\dots,q_{n},p_{n},\dots,p_{_1},
\dots,(p_{_1})^{^{(m-2n+1)}}).\cr}
$$
On $\S$ it reduces to $(p_{_1})_x+{\dd { H}\over{\dd q_{_1}}}\equiv
0$, and the system (2.5) is a canonical Hamiltonian system, with Hamiltonian 
function
$ H$, giving the reduced $x$--flow.

Now we will show that also the $t$--flow reduces on $\S$ with Hamiltonian
function $(-\tilde Q)$.
\bigskip
Firstly we observe that $ Q$ is a function of $x$,$t$,
$u$ and its $x$--derivatives up to the order $(m+n)$, then it can be
rewritten in terms of $x$, $t$, $(p_i,q_i)$ and $p_{{_{1}}}^{{_{(l)}}}$ up to
the order $l=m-n+1$.

We denote with $\tilde f$ a function
$f(x,t,u(x),\dots,u(x)^{(j)})$ reduced on $\S$; notice that, if $j\geq 2n$,
then the reduction can be done using the relation  $$
(p_{_1})_x=-{\dd { H^{^{(a)}}}\over{\dd q_{_1}}}.\eqno{(2.6)}
$$
Then $\tilde f$ does depend explicitly only on the $p_i$ and $ q_i$, for
$i=1,\dots, n$ and on the time $t$. In fact differentiating (t.6) one obtains
the derivatives $p_{{_{1}}}^{{_{(l)}}}$ in terms of the canonical coordinates
$(p_i,q_i)$.
\bigskip
\bigskip
 2) We consider the derivative 
$$
\eqalignno{{dL\over dt}=&{\dd L\over \dd t}+\sum_{i=1}^{n}{\dd L\over \dd
q_i} {dq_i\over dt}+\sum_{i=1}^{n}{\dd L\over \dd p_i}
{dp_i\over dt}=\cr
=&-{\dd H\over \dd t}-\sum_{i=1}^{n}{\dd H\over \dd
q_i} {dq_i\over dt}+\sum_{i=1}^{n}p_i
{d^2q_i\over dx dt}.&(2.7)\cr }
$$
 From the fact that $I$ is a first integral, one deduces that
${dL\over dt}$ must be the total derivative in $x$ of  a
functional $\Lambda$ that does depend
  on $x$,$t$,$(p_i,q_i)$ and
 $p_{{_{1}}}^{{_{(l)}}}$ up to the order
$l=m-n+1$; we have:
$$
\eqalignno{{d\Lambda\over dx}=&{\dd Q\over \dd x}+\sum_{i=1}^{n}{\dd
Q\over \dd q_i} {dq_i\over dx}+\sum_{i=1}^{n}{\dd Q\over \dd p_i}
{dp_i\over dx}+\sum_{i=1}^{m-n+1}{\dd Q\over \dd
(p_{_1})^{(i)}} {d\over
dx}(p_{_1})^{(i)}+\cr
+&\sum_{i=1}^{n}{dp_i\over dx}{dq_i\over
dt}+\sum_{i=1}^{n}p_i {d^2q_i\over dx dt}=\cr
=&{\dd Q\over \dd x}+\sum_{i=1}^{n}{\dd
Q\over \dd q_i} {\dd H\over \dd p_i}-\sum_{i=2}^{n}{\dd Q\over \dd
p_i} {\dd H\over \dd q_i}+\sum_{i=1}^{m-n+1}{\dd Q\over \dd
(p_{_1})^{(i)}} {d\over
dx}(p_{_1})^{(i)}+\cr
-&\sum_{i=2}^{n}{\dd H\over \dd q_i}{dq_i\over
dt}+\sum_{i=1}^{n}p_i {d^2q_i\over dx dt}+{dq_1\over
dt}(p_{_1})_x+{\dd Q\over \dd p_1}(p_{_1})_x.&(2.8)\cr
} $$
Then equation (2.2) gives:
$$
\eqalign{{\dd H\over \dd t}+&{\dd Q\over \dd x}+\sum_{i=1}^{n}{\dd
Q\over \dd q_i}{\dd H\over \dd p_i}-\sum_{i=2}^{n}{\dd Q\over \dd
p_i} {\dd H\over \dd q_i}+{\dd Q\over \dd p_1}(p_{_1})_x
+\cr
+&\sum_{i=1}^{m-n+1}{\dd Q\over \dd
(p_{_1})^{(i)}} {d\over
dx}(p_{_1})^{(i)}+{dq_1\over
dt}(p_{_1})_x+{\dd H\over \dd q_1}{dq_1\over
dt}=0,}
$$
which can be rewritten as 
$$
{\dd H\over \dd t}+{dq_1\over
dt}\biggl((p_{_1})_x+{\dd H\over \dd q_1}\biggr)=-{d\over dx} Q\eqno
{(2.9)}
$$
At this point we need the
\bigskip
\noindent{\bf Lemma 2.1}: {\it On the submanifold $\S$ the following relation
holds}: $$
\widetilde{\biggl(
{\dd Q\over\dd (p_{_1})^{(j)}}\biggr)}=0\qquad \forall j\geq1.\eqno{(2.10)}
$$
\bigskip
\noindent{\bf Proof}: See Appendix 2.A
\bigskip
\bigskip

 Hence, on the submanifold $\S$  eq. (2.9) reduces to: 
$$
{\dd H\over \dd t}+{\dd \tilde {Q}~\over \dd x}+\sum_{i=1}^{n}{\dd
\tilde {Q}~\over \dd q_i}{\dd H\over \dd p_i}-\sum_{i=1}^{n}{\dd \tilde
{Q}~\over \dd p_i} {\dd H\over \dd q_i}=0
$$
This is a {\it zero--curvature} equation:
$$
\{(-\tilde{Q}~),  H\}+{\dd (-\tilde {Q}~)
\over \dd x}-{\dd  H\over \dd t}=0.
$$
This completes  the second step in the proof of the theorem.
\bigskip
{\it 3.}\quad Now we will construct the Hamiltonian system inductively; to this
end we need a further lemma:
\bigskip
\noindent{\bf Lemma 2.2}:  {\it The  fundamental relation} 
$$
{d q_{_1}\over dt}= -{\dd \tilde{Q}~\over {\dd p_{_1}}}.\eqno{(2.11)}
$$
{\it holds}.
\bigskip
\noindent{\bf Proof}: See Appendix 2.A
\bigskip

For simplicity, here and in the following we omit the {\it``tilde"} sign:
$Q$ will indicate the reduced function  on $\S$.
 Now, we assume that ${d q_{_i}\over
dt}= -{\dd Q\over {\dd p_{_i}}}
=-\{q_{_i},Q\}$ and we prove inductively that the same relation holds for
$q_{_{i+1}}$. The scheme of the procedure is the same as in [BN], the
only differences are the contributions of the partial derivatives in $t$ and
$x$.  Indeed, $$ {d q_{_{i+1}}\over dt}=({d q_{_i}\over dt})_x=-
{d\over dx}\{q_{_i},Q\}=-\{\{q_{_i},Q\},H\}-\{q_{_i},
{\dd Q
\over \dd x}\}.
$$
Using the Jacobi identity and the zero--curvature relation we get
$$
{d q_{_{i+1}}\over dt}=-\{{\dd H
\over \dd t}, q_{_i}\}-\{q_{_{i+1}},Q\},\qquad i=1,2,...,n-1.
$$
Here the term $\{{\dd H
\over \dd t}, q_{_i}\}$ is zero for every $i\not= n$ since 
$
\widetilde{\bigl({\dd  L\over \dd t}\bigr)}=-{\dd 
H\over \dd t}
$. Indeed $ L$ depends on $u$ and on the derivatives of $u$ up to 
the order $n$.
This means that,  restricted on $\S$, it  depends on $q_{_1},q_{_2},\dots,
q_{n+1}$. Then, there is no dependence on the $p_{_{i}}$ for $i\not= n$.
Finally we get
$$
\eqalign{\{{\dd H
\over \dd t}, q_{n}\}=&
\{\widetilde{\bigl({\dd  L\over \dd t}\bigr)}, q_{n}\}\not=0,\cr
\{{\dd H
\over \dd t}, q_{i}\}=&0 \qquad i<n.\cr}
$$ 
\bigskip
Hence we have proved that

$$
{d q_{_{i}}\over dt}=-\{q_{_{i}},Q\},\qquad i=1,2,....,n.
$$
\bigskip
\noindent Now we prove that ${d p_{_i}\over dt}=-\{ p_{_i},Q\}$ by
induction, starting from $p_{n}$.

This comes from the commutativity of the flows 
$$ {d \over dt}{d \over dx}q_{n}=
{d \over dx}{d \over dt}q_{n};
$$ 
explicitly:

$$
\eqalign{
{d \over dt}\biggl({d \over dx}q_{n}\biggr)=&
{d \over dt}\biggl({\partial H\over\partial  p_{n}}\biggr)=\cr
=&{\partial\over \partial t}\ {\partial H\over\partial  p_{n}}-
\sum _{i=1}^{n}{\partial H\over\partial  p_{n}\partial q_{_i}}
{\partial Q\over\partial  p_{_i}}+
{\partial^2 H\over\partial  p^2_{n}}{d \over dt}p_{n}=\cr
=&\{q_{n},{\partial H\over\partial t}\}-
\sum _{i=1}^{n}{\partial H\over\partial  p_{n}\partial q_{_i}}
{\partial Q\over\partial  p_{_i}}+
{\partial^2 H\over\partial  p^2_{n}}{d \over dt}p_{n}.\cr}
$$
On the other hand, using the Jacobi identity and the {\it zero--curvature}
equation, one can write 
$$
\eqalign{
{d \over dx}\biggl({d q_{n}\over dt}\biggr)=&
-{\partial\over \partial x}\ {d Q\over dt}-
\{\{q_{n},Q\},H\}=\cr
=&-\{q_{n},\{Q,H\}\}-\{\{q_{n},H\}Q\}-\{q_{n},{\partial Q
\over \partial x}\}=\cr
=&\{q_{n},{\partial H\over\partial t}\}-
\{{\partial H\over \partial p_{n}},Q\}=\cr
=&\{q_{n},{\partial H\over\partial t}\}-
\sum _{i=1}^{n}{\partial H\over\partial  p_{n}\partial q_{_i}}
{\partial Q\over\partial  p_{_i}}+
{\partial^2 H\over\partial  p^2_{n}}{\partial Q\over\partial 
q_{n}}.\cr}
 $$
Comparing the two expressions and noticing that 
${\partial^2 H\over\partial  p^2_{n}}\not=0$ because of the nondegeneracy, we
get
 $$
{d p_{n}\over dt}=-\{ p_{n},Q\}={\partial Q\over\partial 
q_{n}}.
$$-
Now we suppose that ${d p_{_i}\over dt}=\{ p_{_i},Q\}$ and we deduce
the same for $p_{_{i-1}}$. Indeed,
$$
\eqalign{{d \over dx}\biggl({d p_{_i}\over dt}\biggr)
=&\{{\partial Q\over
\partial q_{_i}},H\}+{\partial\over \partial x}{
\partial Q\over \partial
q_{_i}}=\cr
=&-\{\{p_{_i},Q\},H\}-\{p_{_i},{\partial Q\over \partial
x}\}=\cr
=&-\{\{p_{_i},H\},Q\}+\{p_{_i},{\partial H\over \partial
t_{_k}}\}=\cr
=&\{p_{_i},{\partial H\over \partial
t_{_k}}\}+
\sum _{j=1}^{n-1}{\partial^2 H\over\partial  q_{_i}\partial q_{_j}}
{\partial Q\over\partial  p_{_j}}-
{\partial^2 H\over\partial  q_{_i}\partial p_{n}}{\partial Q\over\partial 
q_{n}}-{\partial^2 H\over\partial  q_{_i}\partial p_{_{i-1}}}
{\partial Q\over\partial 
q_{_{i-1}}}\cr}
 $$
and
$$
\eqalign{{d \over dt}\biggl({d p_{_i}\over dx}\biggr)=&
-{d \over dt}\biggl({\partial H\over \partial q_{_i}}\biggr)=\cr
=&
-{\partial\over
\partial t_{_k}}{\partial H\over \partial q_{_i}}-\sum _{j=1}^{n-1}{\partial^2
H\over\partial  q_{_i}\partial q_{_j}} {d q_{_j}\over dt_{_k}}-
{\partial^2 H\over\partial  q_{_i}\partial p_{n}}
{d p_{n}\over dt_{_k}}-{\partial^2 H\over\partial 
q_{_i}\partial p_{_{i-1}}} {d p_{_{i-1}}\over dt_{_k}}=\cr
 =&
\{p_{_i},{\partial H\over \partial
t_{_k}}\}+\sum _{j=1}^{n-1}{\partial^2 H\over\partial  q_{_i}\partial q_{_j}}
{\partial Q\over\partial  p_{_j}}-
{\partial^2 H\over\partial  q_{_i}\partial p_{n}}
{\partial Q\over\partial 
q_{n}}-{\partial^2 H\over\partial 
q_{_i}\partial p_{_{i-1}}} {d p_{_{i-1}}\over dt_{_k}}, 
\cr} 
$$
where ${\partial^2 H\over\partial 
q_{_i}\partial p_{_{i-1}}}=1$.
Comparing the two expressions we get
$
{d p_{_{i-1}}\over dt}={\partial Q\over\partial 
q_{_{i-1}}};
$
hence it follows that
 $$
{d p_{_{i}}\over dt}={\partial Q\over\partial 
q_{_{i}}}=-\{ p_{_{i}},Q \},\qquad i=1,2,...,n
$$
\hfill{Q.E.D.}
\bigskip
\bigskip
\noindent {\bf Remark}: The definition of $Q$:
$$
-Q=-\Lambda +\sum ^{n}_{i=1}p_i {d q_i\over dt}
$$
 looks very similar to the
definition of the Hamiltonian function $H$ of the $x$--flow:
$$
H=-L +\sum ^{n}_{i=1}p_i {d q_i\over dx}
$$
\noindent Here a symmetry between  $x$ and $t$ seems to appear: one could be tempted
to read the definition of $Q$ as a Legendre transform and hence to read 
 $\Lambda$ as the Lagrangian of the $t$--flow. But it is not completely true:
indeed the coordinates $q_i$ and $p_i$ are obtained from the Lagrangian $L$,
they are not, a priori, good coordinates for $\Lambda$. In the next chapter we
will perform a change of coordinates on $\S$, in order to read $\Lambda$ as 
Lagrangian function.

\bigskip
\noindent{\scap 2.A  Appendix}
\bigskip
\noindent{\bf Proof of Lemma 2.1}:
We observe that the recursive relation
$$
\biggl(
{\dd Q\over\dd (p_{_1})^{(j-1)}}\biggr)={\dd\over\dd (p_{_1})^{(j)}}\biggl(
{d\over dx} Q\biggr)-{d\over dx}\biggl(
{\dd Q\over\dd (p_{_1})^{(j)}}\biggr)
$$
holds for $j>1$. Indeed
$$
\eqalign{{\dd\over\dd (p_{_1})^{(j)}}\biggl(
{d\over dx}
Q\biggr)&=\sum_{i=1}^{n}\biggl({\dd^2Q\over\dd q_i\dd 
(p_{_1})^{(j)}}\biggr)(q_i)_x+\sum_{i=1}^{n}\biggl({\dd^2Q\over\dd
p_i\dd  (p_{_1})^{(j)}}\biggr)(p_i)_x+\cr
 &+\sum_{i=1}^{m-n+1}\biggl({\dd^2Q\over\dd
(p_1)^{(i)}\dd  (p_{_1})^{(j)}}\biggr)(p_1)^{(i)}_x+
\biggl({\dd Q\over\dd
  (p_{_1})^{(j-1)}}\biggr).\cr}$$
When we reduce on $\S$:
$$
\widetilde{\biggl(
{\dd Q\over\dd (p_{_1})^{(j-1)}}\biggr)}=-\widetilde{{d\over dx}\biggl(
{\dd Q\over\dd (p_{_1})^{(j)}}\biggr)}.
$$
But $Q$ depends on $(p_{_1})^{(j)}$ up to a finite order, then
$$
\widetilde{\biggl(
{\dd Q\over\dd (p_{_1})^{(j)}}\biggr)}=0\qquad \forall j\geq1.
$$
\hfill{Q.E.D.}
\bigskip
\bigskip
\bigskip
\noindent{\bf Proof of Lemma 2.2}: The expansion of  ${d\over dx} Q$ in
powers of $(p_{_1})_x$ near the point ${\dd  H\over{\dd q_{_1}}}$ reads
$$
\widetilde{\biggl({d\over dx} Q\biggr)}+\widetilde{\biggl[{d\over d(p_{_1})_x}
\biggl({d\over dx} Q\biggr)\biggr]} \biggl((p_{_1})_x+{\dd H\over \dd
q_1}\biggr)+\Theta \biggl((p_{_1})_x+{\dd H\over \dd q_1}\biggr)^2.
 $$
 The zero order term  is
$$
\widetilde{\biggl({d\over dx} Q\biggr)}=-{\dd H\over \dd t}
$$
 by virtue of the zero--curvature equation.
The first order coefficient is
$$
\eqalignno{\widetilde{\biggl[{d\over d(p_{_1})_x}
\biggl({d\over dx} Q\biggr)\biggr]}&=\widetilde{\biggl[{d\over dx}\biggl(
{\dd\over\dd (p_{_1})_x}Q\biggr)\biggr]}+\widetilde{\biggl(
{\dd Q\over\dd p_{_1}}\biggr)}\cr
&=\widetilde{\biggl[{d\over dx}\biggl(
{\dd\over\dd (p_{_1})_x}Q\biggr)\biggr]}+
{\dd\tilde Q\over\dd p_{_1}}-
\sum_{i=2}^{m-n+1}\widetilde{\biggl({\dd Q\over \dd
(p_{_1})^{(i)}} {\dd (p_{_1})^{(i)}\over
\dd p_1}\biggr)},\cr}
 $$

where the only non zero term is ${\dd\tilde Q\over\dd p_{_1}}$, by virtue of
Lemma 2.1. Hence we obtain the power series expansion of ${d\over dx}
Q$ up to the first order:
$$
-{\dd H\over \dd t}+{\dd\tilde Q\over\dd p_{_1}}
 \biggl((p_{_1})_x+{\dd H\over \dd
q_1}\biggr)
$$
this, compared with the left--hand side of eq. (2.9), gives the relation (2.11).

\hfill{Q.E.D.}
\bigskip 

\bigskip
\noindent {\scap 3. Scaling reductions of evolutionary systems:}

\noindent {\scap Lagrangian formulation}
\bigskip
\bigskip
\bigskip
\noindent{\sscap 3.1 General framework}
\bigskip
\noindent
The basic idea is to develop a reduction method dealing on the same
footing  with $x$ and $t$.
The starting point is always the evolutionary PDE
$$
u_t=F(u,u_x,.....,u^{(m)}),\eqno{(3.1)}
$$ 
in the space $\M$ described in Section 1.1.
 The first step of our construction is to read $u$ as a function of
$x$ and $t$ and to consider equation (3.1) as a definition of $u^{(m)}(x,t)$ in
terms of $u(x,t),u_x(x,t),.....,u^{(m-1)} (x,t)$ and $u_t(x,t)$.

This corresponds to consider as  ``coordinates" in $\M$, instead of $u(x,t)$ and
its derivatives in $x$:
$$
u, u_x, u_{xx}, \ldots
$$ 
(here and in the following  $u$ indicate the function $u(x,t)$), the system  
$$
u,u_x,\ldots,u^{(m-1)}, u_t,u_{xt},\ldots,
u^{(m-1)}_t,u_{tt},\ldots
$$ 
By virtue of the reversibility of (3.1) in $u^{(m)}(x,t)$ it is possible to
perform this ``change of variables".
If one introduce the vector
$$
\bar u=(u,u_x,\ldots,u^{(m-1)}),
$$
the new system of ``coordinates" in $\M$ is given by $\bar u(x,t)$ and its 
derivatives in $t$:
$$
\bar u,\bar u_t,\bar u_{tt},\ldots
$$
\bigskip
\noindent At this point one takes a first integral of eq (3.1), i.e. a
functional  
$$
I=\int L\biggl(x,t,u(x),u_x(x),\ldots,u^{(n)}(x)\biggr)dx\eqno{(3.2)}
$$
in the space $\M$, such that
$$
\var I {u(x)}=0.\eqno{(3.3)}
$$
This Euler--Lagrange equation defines a finite dimensional manifold $\S$, i.e.
the set of the fixed points of $I$. Indeed the
Euler--Lagrange equation (3.3) is an ODE of order $2n$, so that the space of the
solutions is a $2n$--dimensional manifold; $\S$ is modeled on this space,
having as coordinates certain combinations of the initial values, i.e.  of
the first $(2n-1)$ $x$--derivatives of $u(x)$ evaluated at $x_0$.
\bigskip
\noindent In Lemma 3.2 below, we rewrite the definition of the manifold $\S$ in
terms of $\bar u(t),\bar u_t(t),\ldots$ and of a functional
$$
J=\int\Lambda\biggl(x,t,\bar u(t),\bar u_t(t),\ldots,\bar u^{(\beta)}(t)\biggr)dt,
$$
where $\Lambda$ can be calculated from $L$ (see eq. (3.4)), and the order 
$\beta$
of derivation in $t$ depends on the ratio between $m$ and $n$, 
as we will show in
detail in Section 3.3.

\noindent In Theorem 3.1 we will prove that $\Lambda\biggl(x,t,\bar u(t),\bar
u_t(t),\ldots, \bar u^{(\beta)}(t)\biggr)$ is the generalized Lagrangian for
the $t$--flow reduced on $\S$. Indeed equation (3.1) can be rewritten in form
of a Euler--Lagrange equation:
$$
\var {J} {\bar u(t)}=0,
$$
for the vector
$$
\var {J} {\bar u(t)}=\biggl(\var {J} {u(t)},\var {J} {u_x(t)},\ldots,
\var {J} {u^{(m-1)}(t)}\biggr),
$$
where
$$
\var {J} {u^{(i)}(t)}={\partial \Lambda\over
\partial u^{(i,0)}} +\sum (-1)^{^\alpha}
{d^{^\alpha}\over {dt^{^\alpha}}}
 {\partial \Lambda\over {\partial {u^{(i,\alpha)}}}}.
$$
In the multiindex $(i,\alpha)$ the Latin character indicates the order in the 
$x$--derivative, the Greek indicate the order in the $t$--derivative.

\noindent Explicitly, equation (3.1) reads
$$
\eqalign{{\partial \Lambda\over
\partial u} +\sum (-1)^{^\alpha}
{d^{^\alpha}\over {dt^{^\alpha}}}
 {\partial \Lambda\over {\partial {u^{(\alpha)}}}}&=0\cr
{\partial \Lambda\over
\partial u_x} +\sum (-1)^{^\alpha}
{d^{^\alpha}\over {dt^{^\alpha}}}
 {\partial \Lambda\over {\partial {u_x^{(\alpha)}}}}&=0\cr
\vdots&\cr
{\partial \Lambda\over
\partial u^{(m-1)}} +\sum (-1)^{^\alpha}
{d^{^\alpha}\over {dt^{^\alpha}}}
 {\partial \Lambda\over {\partial {u^{(m-1,\alpha)}}}}&=0.\cr}
$$
\bigskip
\bigskip
\noindent
We will formalize these facts in the following Theorem 3.1; here we give an
idea of the proof, ignoring all the  calculations, that we
will concentrate in Lemma 3.1 and 3.2, the proof of which is postponed in
Appendix 3.A.

The proof is performed in four steps: firstly we define the new system of 
``coordinates" in $\M$, and we give some useful relations between the new and
the old   ``coordinates". As a second step we rewrite the Lagrangian density
$L\bigl(u(x)\bigr)$ in terms of the new ``coordinates" and we construct, starting from
$\hat L\bigl(\bar u(t)\bigr)$, the Lagrangian density $\hat \Lambda\bigl(\bar
 u(t)\bigr)$.

The third step consists in recovering the relation between  performing the
variation of  
 $L$ (in $x$) and  of  $\hat\Lambda$ (in $t$). The most relevant relation is
that given in Lemma 3.1. This relation is necessary to rewrite the
Euler--Lagrange equation (3.2), defining $\S$, as a condition
on $\hat\Lambda$. The explicit form of this condition is given in Lemma 3.2.

Finally we prove that, under this condition, i.e. after
performing the reduction on $\S$, the starting evolution equation (3.1), reads  
as an Euler--Lagrange equation for $\hat\Lambda$.
\bigskip
\noindent The method of Hamiltonian reduction described in Chapter 2 
allows us to put a canonical system of coordinates $\{p_i,q_i\}$ on $\S$ (see
formula (2.4)). These coordinates are obtained from $L$ via generalized
Lagrange transform, so that they are, in a certain sense, adapted to the
$x$--flow. This means that in these coordinates the reduced $x$--flow is a
Hamiltonian system. Theorem 2.1 also gives the explicit form of the Hamiltonian
function $$
H=-L+\sum_i p_i(q_i)_x.
$$  
The method of Lagrangian reduction which we describe in this Chapter, still
allows us to define a  system of canonical coordinates: we will call it
$\{\tilde p_i,\tilde q_i\}$. These coordinates are obtained from
$\hat\Lambda$, i.e. they are   adapted to the $t$--flow; in fact we will prove
(see Section 3.3) that, in these coordinates, the reduced $t$--flow is a
Hamiltonian system, with Hamiltonian function $$ -\hat Q=-\hat
\Lambda+\sum_i\tilde p_i(\tilde q_i)_t. $$
When rewritten in terms of $\{p_i,q_i\}$, the Hamiltonian $\hat Q$ coincides
with the Hamiltonian function $Q$ constructed by Bogoyavlenskii and Novikov.
\bigskip
\noindent In this sense the alternative definition of $Q$ given by us in 
Theorem 2.1:
$$
-Q=-\Lambda+\sum_i p_i(q_i)_t,
$$ 
is a Legendre transformation, if one uses  the right
system of canonical coordinates (see below).
\bigskip
\bigskip
\noindent{\scap 3.2 Lagrangian reduction}
\bigskip

\noindent {\bf Theorem 3.1:} {\it If the evolutionary PDE:
$$
u_t=F(u,u_x,.....,u^{(m)}),\eqno{(3.3)}
$$
admits a   nondegenerate {\it scaling} symmetry,
then,
 on the manifold $\S$ of the stationary points of the symmetry:
$$
\var I {u(x)}=0,
$$

$$
I=\int L(x,t,u,u_x,\ldots,u^{^{(n)}})
 dx,\qquad {dI\over dt}\equiv 0,
$$
 it reduces to  a
Lagrangian motion in $t$, for the
time dependent Lagrangian function $\Lambda$, determined by:}
$$
{dL\over dt}={d\Lambda\over dx}.\eqno{(3.4)}
$$
 
\bigskip
\noindent {\bf Proof:} (in the following we consider  the case $m\leq n<2m$. The same holds in the case $(\alpha-1)
m\leq n<\alpha m$, as we will show in Section 3.3). We prove the theorem in four
steps:

\bigskip
\noindent{\bf 1.}\quad{\it Change of ``coordinates"}: Let us assume that the
evolutionary equation (3.3) depends on $u(x)$ and on its $x$--derivatives up to
finite order $m$, and that this equation is invertible in $u^{(m)}$. In this
case we can read (3.3) as a definition of $u^{(m)}$ in terms of $
u, u_x,
\dots, u^{(m-1)}$ and $u_t$:
$$
u^{(m)}=f_0\biggl(u, u_x, \dots, u^{(m-1)},u_t\biggr).\eqno{(3.5)}
$$

Differentiating eq. (3.5)  in $x$ one obtains all the  $x$--derivatives of $u$
of order greater then $m$ in terms of $
u, u_x,
\dots, u^{(m-1)}$ and their $t$--derivatives:
 $$
\cases{u^{(m+1)}=
f_{1}\bigl(u, u_x, \dots,
u^{(m-1)},u_t,u_{xt}\bigr)\cr
\vdots&\cr
u^{(2m-1)}=f_{m-1}\bigl(u, u_x, \dots,
u^{(m-1)},u_t,u_{xt}\ldots, u^{(m-1)}_t\bigr)\cr 
u^{(2m)}=f_{m}\bigl(u, u_x, \dots,
u^{(m-1)},u_t,u_{xt}\ldots, u^{(m-1)}_t,u_{tt}\bigr)\cr
\vdots&\cr
u^{(m+n)}=f_{n}\bigl(u, u_x, \dots,
u^{(m-1)},u_t,u_{xt}\ldots,
u^{(m-1)}_t,u_{tt},u_{xtt},\ldots,u^{(n-m)}_{tt}\bigr).
 \cr} 
$$
Explicitly, the first relation has the form
$$
u^{(m+1)}={\partial u^{(m)}\over \partial x}+
{\partial u^{(m)}\over \partial u}u_x+ \dots+
{\partial u^{(m)}\over \partial u^{(m-1)}}u^{(m)}+
{\partial u^{(m)}\over \partial u_t}u_{xt},
$$
and in general,  using the multiindex notation introduced in 
Section 3.1, 
$$
u^{(m+j)}={\partial u^{(m+j-1)}\over \partial x}+\sum_{k=0}^{m-1}
\sum_{\beta=0}^\alpha
{\partial u^{(m+j-1)}\over \partial u^{(k,\beta)}}u^{(k+1,\beta)},
$$
where the higher order $\alpha$ in the $t$--derivative is fixed by 
$(\alpha-1) m\leq j< \alpha m$.
\bigskip 
\noindent This completes the
construction of the map from $$ u, u_x,
\dots, u^{(m-1)},u^{(m)},\ldots,u^{(n)},\ldots,
u^{(2m-1)},u^{(2m)},\ldots\ldots  
$$ 
to the new system of ``coordinates"
$$
u, u_x, \dots,
u^{(m-1)},u_{t},\ldots,u^{(m-1)}_t,u_{tt},u_{xtt},\ldots,u^{(n-m)}_{tt},\ldots
\ldots 
$$
Here below we list some noteworthy relationships between the two system (they
will be useful in the following):
$$
u^{(m)}_t={\partial u^{(m)}\over \partial t}+
{\partial u^{(m)}\over \partial u}u_t+ \dots+
{\partial u^{(m)}\over \partial u^{(m-1)}}u^{(m-1)}_t+
{\partial u^{(m)}\over \partial u_t}u_{tt},\eqno{(3.6a)}
$$

$$
\eqalignno{{d\over dx}\biggl({\partial u^{(i)}\over \partial u_t}\biggr)&=
{\partial u^{(i+1)}\over \partial u_t}-{\partial u^{(i)}\over \partial
u^{(m-1)}}{\partial u^{(m)}\over \partial u_t}&(3.6b)\cr
{d\over dx}\biggl({\partial u^{(i)}\over \partial u^{(k)}_t}\biggr)&=
{\partial u^{(i+1)}\over \partial u^{(k)}_t}-{\partial u^{(i)}\over \partial
u^{(k-1)}_t}&(3.6c)\cr
{d\over dx}\biggl({\partial u^{(i)}\over \partial u^{(k)}}\biggr)&=
{\partial u^{(i+1)}\over \partial u^{(k)}}-
{\partial u^{(i)}\over \partial u^{(k-1)}}-
{\partial u^{(i)}\over \partial
u^{(m-1)}}{\partial u^{(m)}\over \partial u^{(k)}}&(3.6d)\cr
{d\over dt}\biggl({\partial u^{(i)}\over \partial u_t}\biggr)&=
{\partial u^{(i)}_t\over \partial u_t}-{\partial u^{(i)}\over \partial
u}.&(3.6e)\cr}
$$

\bigskip
\bigskip
\noindent{\bf 2.}\quad{\it  Lagrangian densities}:
 The Lagrangian $L$ defining the symmetry, depends on $
u, u_x, \dots, 
u^{(n)}$, so that its derivative ${dL\over dt}$ depends on 
 $
u, u_x, \dots,
u^{(m+n)}$. In terms of the new
``coordinates" one may rewrite $L$ as
$$
\hat L (x,t,u, u_x, \dots,u^{(m-1)},
u_t,\ldots,u^{(n-m)}_t)\eqno{(3.7)}
$$ 
and 
$$
\eqalignno{\widehat{\biggl({dL\over dt}\biggr)}=&\widehat{\biggl(
{\partial L\over \partial t}\biggr)}+\widehat{\biggl({\partial
L\over \partial u}\biggr)} u_t+\dots+
\widehat{\biggl({\partial L\over \partial u^{(m-1)}}\biggr)}
u^{(m-1)}_t+\cr
+&\widehat{\biggl({\partial L\over \partial u^{(m)}}\biggr)}
u^{(m)}_t(u,\ldots,u_{tt})+\widehat
{\biggl({\partial L\over \partial u^{(m+1)}}\biggr)}
u^{(m+1)}_t(u,\ldots,u_{tt},u_{xtt})+\ldots+\cr
+&\widehat{\biggl({\partial L\over \partial u^{(n)}}\biggr)}
u^{(n)}_t(u,\ldots,u_{tt},\ldots,u_{tt}^{(n-m)})&(3.8a).\cr}
$$
Of course, (3.8$a$) coincides with
$$
{d\hat L\over dt}={\partial \hat L\over \partial t}+
{\partial\hat L\over \partial
u} u_t+\dots+{\partial\hat L\over \partial u^{(m-1)}}
u^{(m-1)}_t
+{\partial\hat L\over \partial u_t}
u_{tt}+\ldots+{\partial\hat L\over \partial u^{(n-m)}_t}
u^{(n-m)}_{tt},\eqno{(3.8b)}
$$
where
$$
\eqalignno{{\partial \hat L\over \partial t}=&\widehat{\biggl(
{\partial L\over \partial t}\biggr)}+
\sum_{k=m}^n
\widehat{\biggl(
{\partial L\over \partial u^{(k)}}\biggr)}{\dd u^{(k)}\over\dd
t}&(3.9a)\cr
{\partial\hat L\over \partial
u^{(i)}}=&\widehat{\biggl({\partial L\over \partial
u^{(i)}}\biggr)}+ \sum_{k=m}^n \widehat{\biggl(
{\partial L\over \partial u^{(k)}}\biggr)}{\dd u^{(k)}\over\dd
u^{(i)}}\qquad i=0,\ldots,m-1&(3.9b)\cr
{\partial\hat L\over \partial
u^{(i)}_t}=&
\sum_{k=m+i}^n \widehat{\biggl(
{\partial L\over \partial u^{(k)}}\biggr)}{\dd u^{(k)}\over\dd
u^{(i)}_t}\qquad i=0,\ldots,m-1.&(3.9c)\cr}
$$
From the fact that $I$ is a first integral, it follows that 
there exists a functional
 $$
\hat \Lambda(x,t,u, u_x, \dots, u^{(m-1)},u_{t},
\ldots,u^{(m-1)}_t,u_{tt},\ldots,u^{(n-m-1)}_{tt}),
$$ 
such that
$$
{d\hat L\over dt}={d\hat \Lambda\over dx},
$$
where
$$
\eqalignno{{d\hat \Lambda\over dx}&={\partial\hat \Lambda\over
\partial x}+{\partial\hat \Lambda\over \partial u} u_x+\dots+{\partial\hat
\Lambda\over \partial u^{(m-2)}} u^{(m-1)}+\cr
&+{\partial\hat
\Lambda\over \partial u^{(m-1)}}
u^{(m)}(u,u_x,\ldots,u^{(m-1)},u_{tt})+{\partial\hat \Lambda\over \partial u_t}
u_{xt}+\dots+{\partial\hat
\Lambda\over \partial u^{(m-2)}_t} u^{(m-1)}_t+\cr
&+{\partial\hat
\Lambda\over \partial u^{(m-1)}_t} u^{(m)}_t(u,\ldots,u_{tt})+
{\partial\hat \Lambda\over \partial u_{tt}}
u_{xtt}+\ldots+{\partial\hat \Lambda\over \partial u^{(n-m-1)}_{tt}}
u^{(n-m)}_{tt}.&(3.10)\cr} 
$$
\bigskip
\bigskip
\noindent{\bf 3.}\quad {\it Variations}: Our aim is to reduce equation (3.3)
on the space $\S$ of the stationary points of $I=\int L dx$. This
finite--dimensional manifold is defined by the Euler--Lagrange
equation 
$$
\var I {u(x)}=0. 
$$
This is a variational equation in the old ``coordinates" $u,u_x,\ldots$; how can
we define the same manifold $\S$ in terms of the new ``coordinates"? We must
express $
\var I {u(x)} 
$ in terms of $\Lambda$ and its variations. To this end we recall that
$$
\var I {u(x)}={\partial L\over
\partial u} +\sum (-1)^j
{d^j\over {dt^j}}
 {\partial L\over \partial {u^{(j)}}}. 
$$
and we first
express the terms 
$$
\var L {u^{(j)}(x)}
$$ 
for $j>0$
in terms of $\hat \Lambda$ and the the new ``coordinates", namely:

\bigskip 
\noindent{\bf Lemma
3.1}: {\it The following recurrence relation holds}: $$
{\partial\hat
\Lambda\over\partial u^{(i)}_t}-{d\over
dt}\biggl({\partial\hat
\Lambda\over\partial u^{(i)}_{tt}}\biggr)=\widehat{\biggl(\var {I}
{u^{(i+1)}(x)}\biggr)}+ \sum_{j=m+1}^n\widehat{\biggl(\var {I}
{u^{(j)}(x)}\biggr)}{\partial u^{(j-i)}\over\partial u^{(i)}}\quad i=
1,\ldots,m-1\eqno{(3.11)}
$$
\bigskip
\noindent{\bf Proof}: see Appendix 3.A
\bigskip
\bigskip
\noindent The proof of  Lemma 3.1 is based on the comparison of (3.8) and (3.10)
and their partial derivatives w.r.t. $u^{(j)}_t$ and $u^{(j)}_{tt}$. With a
similar technique, and using equation (3.11), one can  proves the fundamental 
\bigskip 
\noindent{\bf Lemma 3.2}: {\it The
(generalized) Euler--Lagrange equation  
$$
\var I {u(x)}=0
$$
is equivalent to the condition}:
$$
{\partial\hat \Lambda\over \partial u^{(m-1)}}-{d\over dt}\biggl(
{\partial\hat \Lambda\over \partial u^{(m-1)}_{t}}\biggr)=0.\eqno{(3.12)}
$$
\bigskip
\noindent{\bf Proof}: see Appendix 3.A
\bigskip
\noindent Introducing the functional
$$
J=\int \hat\Lambda
\biggl(x,t,u(t),u_x(t),\ldots,u^{(m-1)}(t),u_t(t),\ldots,u^{(n-m-1)}_{tt}(t)
\biggr) dt,
$$ 
equation (3.12) reads
$$
\var J {u^{(m-1)}(t)}=0.
$$
Notice that the object in the left hand side is the last component of the
vector 
$$
\var J {\bar u(t)}=\biggl(\var J {u(t)},\var J {u_x(t)},\ldots,
\var J {u^{(m-1)}(t)}\biggr).
$$

\bigskip
\noindent{\bf 4.}\quad{\it  Reduced evolutionary equation}:
\noindent Here we prove that all the components of the vector $\var J {\bar
u(t)}$ are zero. This can be done
 recursively, by mean of 
\bigskip
\noindent{\bf Lemma 3.3}: {\it The following recurrence relation holds}:
$$ 
\var {J} {u^{(i-1)}(t)}=\var{J}
{u^{(m-1)}(t)} {\partial u^{(m)}\over  \partial u^{(i)}}-{d\over dx}\biggl(\var
{J} {u^{(i)}(t)}\biggr) \qquad i=1,\ldots,m-1.\eqno{(3.13a)}
$$ 
\bigskip
\noindent{\bf Proof}: see Appendix 3.A
\bigskip
\noindent Indeed,  Lemma 3.2 states that the  $(m-1)$--th component 
$\var {J} {u^{(m-1)}(t)}$
is zero when reduced on $\S$, hence, by virtue of (3.13$a$), all the components
of  $\var {J} {\bar u(t)}$ vanish on $\S$. 

\noindent This is the  Euler--Lagrange equation for the Lagrangian 
$$
\hat\Lambda
(x,t,u,u_x,\ldots,u^{^{(m-1)}},u_t,
\ldots,u^{(m-1)}_t,u_{tt},\ldots,u^{(n-m)}_{tt}).
$$

\hfill{Q.E.D.}
\bigskip
\bigskip 
\noindent{\bf Remark}: Equation (3.13$a$) can be rewritten as 
$$
\eqalignno{&\biggl[
{\partial\hat \Lambda\over \partial u^{(i-1)}}-{d\over dt}\biggl(
{\partial\hat \Lambda\over \partial u^{(i-1)}_{t}}\biggr)
+{d^2\over dt^2}\biggl(
{\partial\hat \Lambda\over \partial u^{(i-1)}_{tt}}\biggr)\biggr]=\cr
&=\biggl[{\partial\hat \Lambda\over \partial u^{(m-1)}}-{d\over dt}\biggl(
{\partial\hat \Lambda\over \partial u^{(m-1)}_{t}}\biggr)\biggr]
{\partial u^{(m)}\over  \partial u^{(i)}}-{d\over dx}\biggl[
{\partial\hat \Lambda\over \partial u^{(i)}}-{d\over dt}\biggl(
{\partial\hat \Lambda\over \partial u^{(i)}_{t}}\biggr)
+{d^2\over dt^2}\biggl(
{\partial\hat \Lambda\over \partial u^{(i)}_{tt}}\biggr)\biggr].&(3.13b)\cr}
$$
{\it for} $i=1,\ldots,m-1$.
And The Euler--Lagrange equation reads
$$
\cases{{\partial \Lambda\over \partial u^{(i)}}-{d\over dt}
{\partial \Lambda\over \partial u^{(i)}_t}=0\qquad
i=n-m,\ldots,m-1\cr 
{\partial \Lambda\over \partial u^{(i)}}-
{d\over dt}\bigl({\partial \Lambda\over \partial u^{(i)}_t}\bigr)+
{d^2\over dt^2}{\partial \Lambda\over \partial u^{(i)}_{tt}}=0
\qquad i=0,\ldots,n-m-1.\cr}\eqno{(3.14)}
$$
\bigskip
\bigskip

\bigskip
\noindent{\sscap 3.3 Relation with the  Hamiltonian reduction}
\bigskip
\bigskip
Theorem 3.1 provides an alternative definition of the space $\S$, and of the
relative system of canonical coordinates: 
$$
\cases{\tilde q_{_i}=q_{_i}=u^{(i-1)}\qquad i=1,\ldots,m\cr
\tilde q_{_{m+i}}=(q_{_i})_t=u_t^{(i-1)}\qquad i=1,\ldots,n-m\cr 
\tilde p_{_i}=\var{J}{u_t^{(i-1)}(t)}\qquad i=1,\ldots,m\cr
\tilde p_{_{m+i}}={\partial \hat\Lambda\over\partial u_{tt}^{(i-1)}}\qquad
i=1,\ldots,n-m\cr} \eqno{(3.15)}
$$
We consider now the Hamiltonian ($-Q)$ of the reduced $t$--flow, defined in
Theorem 2.1 
$$
Q=\Lambda -\sum^n_{i=1}p_i(q_i)_t=\Lambda -\sum^{n-1}_{i=0}\var I {u^{(i+1)}(t)}
u^{(i)}_t
$$
At the end of the Chapter 2 we noticed how this expression looks very similar to
a Legendre transform, but it is not; here we will show that actually the Legendre
transform of the Lagrangian $\hat \Lambda$ gives the Hamiltonian $\hat Q$,
where $\hat Q$ is written in the coordinate system relative to $\hat \Lambda$.
\bigskip

\noindent Firstly we rewrite $Q$ in the coordinate system (3.15):
$$
\eqalign{-\hat Q&=-\hat \Lambda +\sum^{n-1}_{i=0}\widehat{\biggl(\var I
{u^{(i+1)}(x)}\biggr)}\widehat{\biggl(u^{(i)}_t\biggr)}\cr
&=-\hat \Lambda+\sum^{m-1}_{i=0}\biggl[\widehat{\biggl(\var I
{u^{(i+1)}(x)}\biggr)}+\sum^{n}_{j=m+1}\widehat{\biggl(\var I
{u^{(j)}(x)}\biggr)}{\partial u^{(j-1)}\over
\partial u^{(i)}}\biggr]u^{(i)}_t+\cr
&+\sum^{n-m-1}_{i=0}\sum^{n}_{j=m+i+1}\widehat{\biggl(\var I
{u^{(j)}(x)}\biggr)}{\partial u^{(j-1)}\over
\partial u^{(i)}_t}u^{(i)}_{tt}.
\cr} 
$$
Using  Lemma 3.1 we get
$$
-\hat Q=-\hat \Lambda+\sum^{m-1}_{i=0}\var {J} {u_t^{(i)}}u^{(i)}_t
+\sum^{n-m-1}_{i=0}{\partial \hat\Lambda\over\partial u^{(i)}_{tt}}u^{(i)}_{tt}=
-\hat \Lambda+\sum^n_{i=1}\tilde p_i(\tilde q_i)_t,\eqno{(3.16)}
$$ 
for $\Lambda\bigl(x,t,\tilde q_i,(\tilde q_1)_t,\ldots,(\tilde q_{n-m})_t\bigr)$.
\bigskip
\bigskip
\bigskip
\noindent {\scap 3.4 Concluding remarks}
\bigskip
\noindent The case  considered in Theorem 3.1 is the more general
one. Indeed, for $(\alpha-1)m<n\leq \alpha m$, the Euler--Lagrange equation
defining $\S$: $$
\widehat{\biggl(\var I {u(x)}\biggr)}=0
$$
into the new ``coordinates", is a differential equation
in $u,\ldots,u^{(m-1)},\ldots,
u^{(n-m,\alpha)}$.

\bigskip
\noindent The Lagrangian $L$ transforms into
$$\hat L(x,t,u,u_x,\ldots,u^{(m-1)},u_t,\ldots,u^{(m-1)}_t,u_{tt},\ldots,
u^{(n-m,\alpha-1)})
$$ 
and we can define the new Lagrangian 
$$
\hat \Lambda(x,t,u,u_x,\ldots,u^{(m-1)},u_t,\ldots,u^{(m-1)}_t,u_{tt},\ldots,
u^{(n-m-1,\alpha)}).
$$
The proof of the theorem is the same, one has only to consider the
identities $$
{\partial\over\partial u^{(i,\beta)}}\widehat{\biggl({d L\over dt}\biggr)}=
{\partial \over\partial u^{(i,\beta)}}\biggl({d\hat \Lambda\over dx}\biggr)
$$
for $i=1,\ldots,m-1$ and $\beta=1,\ldots,\alpha$.
\bigskip
\noindent In particular, the manifold $\S$ is defined by
$$
\var {J} {u^{(m-1)}(t)}=0\eqno{(3.17)}
$$
and it naturally carries the canonical system of coordinates
$$
\cases{\hat q_{_{\beta m+i}}=u^{(i-1,\beta)}\qquad i=1,\ldots,m;\quad
\beta=0,\ldots,\alpha-2\cr
\hat q_{_{(\alpha-1)m+i}}=u^{(i-1,\alpha-1)}\qquad i=1,\ldots,n-(\alpha-1)m\cr
 \hat p_{_i}=\var{J}{(\hat q_{_i})_t}\qquad
i=1,\ldots,n\cr}\eqno{(3.18)} $$
\bigskip
\bigskip
\noindent In the following we will consider in detail the case $\alpha=1$, i.e.
$n<m$, which occurs in the applications we are interested in (see next
chapter). In this case $L$ and $\hat L$ coincide and  $$
{dL\over dt}={\partial L\over \partial t}+{\partial L\over \partial u}
u_t+\dots+{\partial L\over \partial u^{(n)}}
u^{(n)}_t.\eqno{(3.19)}
$$
The new Lagrangian is 
$$\hat\Lambda(x,t,u, u_x, \dots, u^{(m-1)},(u)_{t},
\ldots,u^{(n-1)}_t)
$$
 with
$$
{d L\over dt}={d\Lambda\over dx}={\partial \Lambda\over \partial
x}+{\partial \Lambda\over \partial u} u_x+\dots+{\partial \Lambda\over \partial
u^{(m-1)}} u^{(m)}+{\partial
\Lambda\over \partial u_t} u_{xt}+\dots+{\partial \Lambda\over \partial
u^{(n-1)}_t} u^{(n)}_t.\eqno{(3.20)}
$$

\bigskip
\noindent In this case Lemma 3.1 reduces to the following recurrence relation: 
$$
\var I {u^{(i)}(x)}={\partial \Lambda\over\partial u^{(i-1)}_t}
\qquad i=1,\ldots,n\eqno{(3.21)}
$$ 
and the proof is based on the identity
$$
{\partial\over\partial u^{(i)}_t}\biggl({dL\over dt}\biggr)=
{\partial \over\partial u^{(i)}_t}\biggl({d\Lambda\over dx}\biggr)
\qquad i=0,\ldots,n-1,\eqno{(3.22)}
$$
observing that, for $i\geq1$,
$$
{\partial\over\partial u^{(i)}_t}\biggl({dL\over dt}\biggr)=
{\partial L\over\partial u^{(i)}}.
$$
In fact, $L$ does not depend on $u^{(i)}_t$,
and
$$
{\partial\over\partial u^{(i)}_t}\biggl({d\Lambda\over dx}\biggr)=
{d\over dx}\biggl({\partial \Lambda\over\partial u^{(i)}_t}\biggr)+
{\partial \Lambda\over \partial u^{(i-1)}_t}.
$$
In particular
 the first step, $i=n$, follow directly from the fact that
 the only dependence of $u^{(n)}_t$ in both (3.19) and (3.20) is
the one explicitly shown, so that
 $$
{\partial L\over \partial u^{(n)}}=
{\partial \Lambda\over \partial u^{(n-1)}_{t}}.\eqno{(3.23)}
$$
On the other hand, from (3.22) for the index $i=0$, one obtains the fundamental
relation
 $$
\var I {u(x)}={\partial \Lambda\over \partial u^{(m-1)}}{\partial u^{(m)}\over 
\partial u_t},\eqno{(3.24)}
$$
and, since ${\partial u^{(m)}\over 
\partial u_t}$ is always nonzero,  the condition that
defines the submanifold $\S$ is 
$$
{\partial \Lambda\over \partial u^{(m-1)}}=0.
$$ 
The relative system of canonical coordinates is given by:
$$
\cases{\hat q_{_i}=u^{(i-1)},\cr 
\hat p_{_i}={\partial\Lambda\over\partial u_t^{(i-1)}}\cr}
$$
For $i=1,\ldots,n$. 
\bigskip
\bigskip
\noindent
The reduced $t$--flow is Lagrangian, with Lagrangian
$\Lambda$.

\noindent Indeed, from the identity
$$
{\partial\over\partial u^{(i)}}\biggl({dL\over dt}\biggr)=
{\partial\over\partial u^{(i)}}\biggl({d\Lambda\over dx}\biggr)
\qquad i=1,\ldots,m-1,
$$
and using the Lemma,
one obtains on the subspace $\S$:

$$
\cases{{\partial \Lambda\over \partial u^{(i)}}=0\qquad i=n,\ldots,m-1\cr
{\partial \Lambda\over \partial u^{(i)}}-
{d\over dt}\bigl({\partial \Lambda\over \partial u^{(i)}_t}\bigr)=0
\qquad i=0,\ldots,n-1\cr}
$$
which is the Euler--Lagrange equation for the Lagrangian $\Lambda
(x,t,u,u_x,\ldots,u^{^{(m-1)}},u_t,
\ldots,u^{(n-1)}_t)$.
\bigskip
\bigskip
\bigskip
\noindent{\scap 3.5 Example: KdV with $t_7$ fixed}
\bigskip
\noindent We will give below an example of how does Theorem 3.1 works for the
first non trivial case, $n=m$. We study  the Lagrangian reduction of the  KdV
equation    
$$
u_t=6uu_x-u_{xxx}\eqno{(3.25)}
$$
 on the stationary manifold of the $t_7$--flow. 

The 
Lagrangian density of the $t_7$--flow, reduced to the normal form (here I
mean that $L$ does not contains total derivatives), depends on the
$x$--derivatives of $u(x,t)$ up to  order $n=3$, and has the expression
$$
L=7u^5+35u^2u_x^2+7uu_{xx}^2+{1\over 2}(u^{(3)})^2.\eqno{(3.26)}
$$   
The submanifold $\S$ of the stationary points  defined by the Euler--Lagrange
equation for $L$ gives the $n+m=6$ derivative in terms of the first five,
explicitly
$$
u^{(6)}=14uu^{(4)}+28u_xu_{xxx}-70u^2u_{xx}+21u_{xx}^2-70uu_x^2+35u^4.
\eqno{(3.27)}
$$
From the relation
$$
{dL\over dt}={d\Lambda\over dx}
$$
one construct the Lagrangian
$\Lambda(u, u_x,\ldots,u^{(5)})$. By direct calculation
$$
\eqalignno{\Lambda=&-u^{(3)}u^{(5)}+{1\over
2}(u^{(4)})^2-14uu_{xx}u^{(4)}+10u(u^{(3)})^2+14u_xu_{xx}u_{xxx}+\cr
-&70u^2u_xu^{(3)}-(u^{(2)})^3+77u^2(u_{xx})^2+70uu_x^2u_{xx}-35u^4u_{xx}+\cr
-&{35\over 2}u_x^4+280u^3u_x^2+35u^6.&(3.28)
\cr} $$
The evolution equation (3.27) is the definition of $u_{xxx}$ in
terms of $(u, u_x,u_{xx},u_t)$, explicitly: 
$$
u_{xxx}=6uu_x-u_t. 
$$
Differentiating this relation in $x$ one
obtains  $$
\cases{u^{(4)}=6uu_{xx}+6u_x^2-u_{xt}\cr
u^{(5)}=18u_xu_{xx}+36u^2-6uu_t-u_{xxt}\cr
u^{(6)}=18u_{xx}^2+180uu_x^2+36u^2u_{xx}-30u_xu_t-12uu_{xt}+u_{tt}\cr}
$$
which is a map from the ``coordinates"
$$
u, u_x,u_{xx}, u^{(3)},u^{(4)},u^{(5)},u^{(6)},\ldots
$$
into
$$
u, u_x,u_{xx},u_t, u_{xt},u_{xxt},u_{tt},\ldots
$$
\bigskip 

 The Lagrangian $L$ depends on $
u, u_x,u_{xx},u_{xxx}$; in the new ``coordinates"
$$
\hat L(u, u_x,u_{xx},u_t)=7u^5+53u^2u_x^2+7uu_{xx}^2-6uu_xu_t+{1\over 2}u_t^2
$$
Its derivative ${d\hat L\over dt}$ looks like
$$
{d\hat L\over dt}={\partial\hat  L\over \partial t}+{\partial\hat L\over \partial
u} u_t+{\partial\hat L\over \partial u_x}
u_{xt}+{\partial\hat  L\over \partial u_{xx}}
u_{xxt}+
{\partial\hat L\over \partial u_t}
u_{tt}.
$$
And there exist a functional 
$\hat \Lambda$ depending on $u, u_x,u_{xx},u_t,u_{xt},u_{xxt}$, explicitly

$$
\eqalign{\hat\Lambda=&35u^6+4u^3u_x^2+{1\over 2}u_x^4-6u_x^2u_{xt}+{1\over
2}u_{xt}^2+\cr
-&35u^4u_{xx}-2uu_x^2u_{xx}+8uu_{xt}u_{xx}+11u^2u_{xx}^2+\cr
-&u_{xx}^3+6uu_xu_{xxt}+22u^2u_xu_t+4u_xu_{xx}u_t-u_{xxt}u_t
+4uu_t^2\cr}
$$
such that
$$
\eqalign{{d\hat L\over dt}&={d\hat\Lambda\over dx}={\partial\hat \Lambda\over
\partial x}+{\partial\hat \Lambda\over \partial u} u_x+
+{\partial\hat \Lambda\over \partial u_x} u_{xx}+{\partial\hat \Lambda\over
\partial u_{xx}} u^{(3)}+\cr
&+{\partial\hat
\Lambda\over \partial u_t} u_{xt}+{\partial\hat
\Lambda\over \partial u_{xt}} u_{xxt}+{\partial\hat
\Lambda\over \partial u_{xxt}} u^{(3)}_t,\cr}
$$
where
$$
\eqalign{{\partial L\over \partial t}&=0\cr
{\partial L\over \partial u}&=35u^4+106uu_x^2+7u_{xx}^2-6u_xu_t\cr
{\partial L\over \partial u_x}&=106u^2u_x-6uu_t\cr
{\partial L\over \partial u_{xx}}&=14uu_{xx}\cr
{\partial L\over \partial u_t}&=-6uu_x+u_t\cr}
$$
and
$$
\eqalign{{\partial\hat \Lambda\over \partial x}&=0\cr
{\partial\hat\Lambda\over \partial
u}&=210u^5+12u^2u_x^2-140u^3u_{xx}-2u_x^2u_{xx}+8u_{xt}u_{xx}
+22uu_{xx}^2+6u_xu_{xxt}+44uu_xu_t+4u_t^2\cr {\partial
\hat\Lambda\over \partial u_x}&=8u^3u_x+2u_x^3-12u_xu_{xt}-4uu_xu_{xx}
+6uu_{xxt}+22u^2u_t+4u_{xx}u_t\cr 
{\partial\hat \Lambda\over \partial
u_{xx}}&=-35u^4-2uu_x^2+8uu_{xt}+22u^2u_{xx}-3u_{xx}^2+4u_xu_t\cr 
{\partial\hat \Lambda\over \partial
u_t}&=22u^2u_x+4u_xu_{xx} -u_{xxt}+8uu_t\cr
{\partial\hat \Lambda\over \partial u_{xt}}&=-6u_x^2+u_{xt}+8uu_{xx}\cr
{\partial\hat \Lambda\over \partial u_{xxt}}&=6uu_x-u_t\cr}
$$
\bigskip
\bigskip

\noindent Lemma 3.2 states that the condition $\var I u=0$, which
defines the submanifold $\S$, is equivalent to the condition
$$
{\partial\hat \Lambda\over \partial u_{xx}}-{d\over dt}\biggl(
{\partial\hat \Lambda\over \partial u_{xxt}}\biggr)=0,
$$ 
explicitly:
$$
u_{tt}=35u^4+2uu_x^2-22u^2u_{xx}+3u_{xx}^2+2u_xu_t-2uu_{xt}.
$$
The Euler--Lagrange equation for $\hat\Lambda$ reads
$$
\eqalignno{{\partial\hat \Lambda\over \partial u_{xx}}-
{d\over dt}\bigl(
{\partial\hat \Lambda\over \partial u_{xxt}}\bigr)
&=u_{tt}-35u^4-2uu_x^2+22u^2u_{xx}-3u_{xx}^2-2u_xu_t+2uu_{xt}=0\cr
{\partial\hat \Lambda\over \partial u_x}-
{d\over dt}\bigl(
{\partial\hat \Lambda\over \partial u_{xt}}\bigr)&=
8u^3u_x+2u_x^3-4uu_xu_{xx}
-2uu_{xxt}+22u^2u_t-4u_{xx}u_t-u_{xtt}
=0\cr
{\partial\hat \Lambda\over \partial u}-
{d\over dt}\bigl({\partial\hat \Lambda\over \partial u_t}\bigr)&=
210u^5+12u^2u_x^2-140u^3u_{xx}-2u_x^2u_{xx}+4u_{xt}u_{xx}
+22uu_{xx}^2+\cr
&+2u_xu_{xxt}-4u_t^2-22u^2u_{xt}+u_{xxtt}-8uu_{tt}=0.&(3.29)
\cr}
$$
\bigskip
\bigskip
\noindent In this case $L$ is nondegenerate, so that on $\S$ we can define the
system of canonical coordinates

$$
\cases{q_i=u^{(i-1)}\qquad i=1,2,3\cr
p_i={\partial\hat \Lambda\over \partial u^{(i-1)}_t}\qquad i=1,2,3\cr}
$$
which reads
$$
\eqalign{\tilde q_1&=u\cr
\tilde q_2&=u_x\cr
\tilde q_3&=u_{xx}\cr
\tilde p_1&=22u^2u_x+4u_xu_{xx} -u_{xxt}+8uu_t\cr
\tilde p_2&=-6u_x^2+u_{xt}+8uu_{xx}\cr
\tilde p_3&=6uu_x-u_t\cr}
$$
\bigskip
\bigskip
\noindent We will now solve the problem from the Hamiltonian point of
view: starting from $L$ and following Theorem 1.1 one construct the canonical
coordinates $\{p_i,q_i\}$ on $\S$:
$$
\eqalign{ q_1&=u\cr
 q_2&=u_x\cr
 q_3&=u_{xx}\cr
 p_1&=\var {I} {u_x}=70 u^2u_x-14u_xu_{xx}-14uu_{xxx}+u^{(5)}\cr
 p_2&=\var {I} {u_{xx}}=14 uu_{xx}-u^{(4)}\cr
 p_3&=\var {I} {u_{xxx}}=u_{xxx}\cr}
$$
and the Hamiltonian function
$$
Q=\Lambda-\sum_{i=1}^3 p_1(q_i)_t.
$$
By direct calculation one obtains
$$
\eqalign{Q&=35u^6-140u^3u_x^2-{35\over 2}u_x^4-35u^4u_{xx}+70uu_x^2u_{xx}-
7u^2u_{xx}^2-u_{xx}^2+84u^2u_xu_{xxx}+\cr
&-18u_xu_{xx}u_{xxx}-10uu_{xxx}^2
+6u_x^2u^{(4)}+6uu_{xx}u^{(4)}-{1\over 2}(u^{(4)})^2-6uu_xu^{(5)}+
u_{xxx}u^{(5)}\cr},
$$
and in canonical coordinates
$$
\eqalign{Q&=35q_1^6+280 q_1^3q_2^2-{35\over 2}q_2^4-35q_1^4q_3+70
q_1q_2^2q_3-21q_1^2q_3^2+\cr
&-q_3^3-6q_1q_2p_1-6 q_2^2 p_2+8q_1q_3p_2-{1\over 2} p_2^2-
70 q_1^2q_2p_3-4q_2q_3p_3+p_1p_3+4 q_1p_3^3.\cr}
$$
The corresponding Hamiltonian system reads
$$
\eqalign{ \dot q_1&=6q_1q_2-p_3\cr
 \dot q_2&=6q_2^2-8q_1q_3+p_2\cr
 \dot q_3&=70q_1^2q_2+4q_2q_3-p_1-8q_1p_3\cr
 \dot p_1&=210q_1^5+840q_1^2q_2^2-140q_1^3q_3+70q_2^2q_3-
42q_1q_3^2-6q_2p_1+8q_3p_2-140q_1q_2p_3+4p_3^2\cr
 \dot p_2&=560q_1^3q_2-70q_2^3+140q_1q_2q_3-6q_1p_1-
12q_2p_2-70q_1^2p_3-4q_3p_3\cr
 \dot p_3&=-35q_1^4+70q_1q_2^2-42q_1^2q^3-3q_3^2+8q_1p_2-4q_2p_3\cr}.
$$ 
Rewriting this system in coordinates $\{\tilde p_i,\tilde q_i\}$ one obtains
exactly (3.29). 
 \bigskip
\bigskip
\bigskip
\noindent{\scap 3.A Appendix}
\bigskip
\bigskip
\noindent{\bf  Proof of Lemma 3.1}:
 we prove the Lemma in two parts:
\bigskip
$\bullet$ firstly we prove the relation
$$
{\partial\hat
\Lambda\over\partial u^{(i)}_{tt}}=\sum_{j=i+m+1}^n\widehat{\biggl(\var {I}
{u^{(j)}(x)}\biggr)}{\partial
u^{(j-i)}\over\partial u^{(i)}_{t}} \qquad i=0,\ldots,n-m-1.\eqno{(a.1)}
$$
\bigskip
\bigskip
\noindent For convenience we can explicitly rewrite eq. (3.8a), using  (3.9):
$$
\eqalignno{\widehat{\biggl({dL\over dt}\biggr)}=&\biggl[\widehat{\biggl(
{\partial L\over \partial t}\biggr)}+
\sum_{j=m}^n
\widehat{\biggl(
{\partial L\over \partial u^{(j)}}\biggr)}{\dd u^{(j)}\over\dd
t}\biggr]+\cr
+\sum_{i=0}^{m-1}&\biggl[\widehat{\biggl({\partial L\over \partial
u^{(i)}}\biggr)}+ \sum_{j=m}^n
\widehat{\biggl(
{\partial L\over \partial u^{(j)}}\biggr)}{\dd u^{(j)}\over\dd
u^{(i)}}\biggr]u^{(i)}_{t}+\cr
+\sum_{i=0}^{n-m}&\biggl[\sum_{j=m+i}^n
\widehat{\biggl(
{\partial L\over \partial u^{(j)}}\biggr)}{\dd u^{(j)}\over\dd
u^{(i)}_t}\biggr]u^{(i)}_{tt}.&(a.2)\cr}
$$
Notice that the arguments of in the square brackets depend on $u$ and its
$x$--derivatives upon the order $m-1$ and on $u_t$ and its $x$--derivatives
upon the order $n-m$, so that the  dependence on $u^{(i)}_t$ for
$n-m+1\leq i\leq m-1$ and on $u^{(j)}_{tt}$, for every $j$ is only the explicit
one.
Analogously
$$
\eqalignno{{d\hat \Lambda\over dx}=&\biggl[
{\partial\hat \Lambda\over
\partial x}+{\partial\hat
\Lambda\over \partial u^{(m-1)}_t}{\dd u^{(m)}\over\dd
t}\biggr]+\sum_{i=1}^{m}{\partial\hat
\Lambda\over \partial u^{(i-1)}}u^{(i)}+\cr
+&\biggl[{\partial\hat
\Lambda\over \partial u^{(m-1)}_t}{\dd u^{(m)}\over\dd
u}\biggr]u_t+\sum_{i=1}^{m-1}\biggl[
{\partial\hat \Lambda\over \partial u^{(i)}_t}+{\partial\hat
\Lambda\over \partial u^{(m-1)}_t}{\dd u^{(m)}\over\dd
u^{(i)}}\biggr] u^{(i)}_t+\cr
+&{\partial\hat
\Lambda\over \partial u^{(m-1)}_t}{\dd u^{(m)}\over\dd
u_t} u_{tt}
+\sum_{i=1}^{n-m}{\partial\hat
\Lambda\over \partial u^{(i-1)}_tt} u^{(i)}_{tt}.&(a.3)\cr} 
$$  
 The $i$--th step of ($a$.1) is obtained
from the obvious identity
$$
{\partial\over\partial u^{(i)}_{tt}}\widehat{\biggl({d L\over dt}\biggr)}=
{\partial \over\partial u^{(i)}_{tt}}\biggl({d\hat \Lambda\over dx}\biggr)
\qquad i=1,\ldots,n-m.
$$
Indeed, from ($a$.2) and ($a$.3), it follows that
$$
{\partial\over\partial u^{(i)}_{tt}}\widehat{\biggl({d L\over dt}\biggr)}=
\sum^{n}_{j=m+i}\widehat{\biggl({\partial L\over \partial u^{(j)}}\biggr)}
{\partial
u^{(j)}\over\partial u^{(i)}_{t}},\eqno{(a.4)}
$$
and
$$
{\partial\over\partial u^{(i)}_{tt}}\biggl({d\hat \Lambda\over dx}\biggr)=
{d\over dx}\biggl({\partial\hat \Lambda\over\partial u^{(i)}_{tt}}\biggr)+
{\partial \Lambda\over \partial u^{(i-1)}_{tt}}.\eqno{(a.5)}
$$
In particular, at the first step, $i=n-m$ one obtains the basic relation
$$
{\partial\hat \Lambda\over \partial u^{(n-m-1)}_{tt}}=
\widehat{\biggl({\partial L\over \partial u^{(n)}}\biggr)}
{\partial u^{(n)} \over \partial u^{(n-m)}_{t}}\equiv
\widehat{\biggl(\var I {u^{(n)}(x)}\biggr)}
{\partial u^{(n)} \over \partial u^{(n-m)}_{t}}
.\eqno{(a.6)}
$$
Substituting ($a$.6) into the further step of the recurrence, one finds
$$
{\partial \Lambda\over \partial u^{(n-m-2)}_{tt}}=
\widehat{\biggl(\var I {u^{(n-1)}(x)}\biggr)}
{\partial u^{(n-2)} \over \partial u^{(n-m-2)}_{t}}+
\widehat{\biggl(\var I {u^{(n)}(x)}\biggr)}
{\partial u^{(n-1)} \over \partial u^{(n-m-2)}_{t}}
$$
and so on. This gives relation ($a$.2).
\bigskip
\bigskip
$\bullet$ The second step is the proof of the relation
$$ 
{\partial\hat
\Lambda\over\partial u^{(i)}_t}=\widehat{\biggl(\var {I} {u^{(i+1)}(x)}\biggr)}+
\sum_{j=m+1}^n\widehat{\biggl(\var {I} {u^{(j)}(x)}\biggr)}{\partial
u^{(j-i)}\over\partial u^{(i)}} \qquad i=n-m,\ldots,m-1,\eqno{(a.7)}
$$
which is a part of ($a$.1), indeed, for $i\geq n-m$, the partial derivative
${\partial\hat
\Lambda\over\partial u^{(i)}_{tt}}$ vanishes.  
Very much as in the previous case, equation ($a$.7) follows from the identity
$$
{\partial\over\partial u_{tt}}\widehat{\biggl({d L\over dt}\biggr)}=
{\partial\over\partial u_{tt}}\biggl({d\hat \Lambda\over dx}\biggr)
$$
Using  ($a$.2) one can rewrite the left hand side as
$$
{\partial\over\partial u_{tt}}\widehat{\biggl({d L\over dt}\biggr)}=
\sum^n_{j=m}\widehat{\biggl({\partial L\over\partial u^{(j)}}\biggr)}
{\partial u^{(j)} \over \partial u_{t}}.\eqno{(a.8)}
$$
On the other hand, from ($a$.3), one has
$$
{\partial\over\partial u_{tt}}\biggl({d\hat \Lambda\over dx}\biggr)=
{d\over dx}\biggl({\partial\hat \Lambda\over\partial u_{tt}}\biggr)+
{\partial\hat \Lambda\over \partial u^{(m-1)}_t}{\partial u^{(m)}\over 
\partial u_{t}},\eqno{(a.9)}
$$
where
$$
{d\over dx}\biggl({\partial\hat \Lambda\over\partial
u_{tt}}\biggr)={d\over dx}
\biggl(\sum_{j=m+1}^n\widehat{\bigl(\var {I}
{u^{(j)}(x)}\bigr)}{\partial
u^{(j-1)}\over\partial u_{t}}\biggr).
$$
We develop the right hand side, recalling (3.6b):
$$
{d\over dx}\biggl({\partial u^{(i)}\over \partial u_t}\biggr)=
{\partial u^{(i+1)}\over \partial u_t}-{\partial u^{(i)}\over \partial
u^{(m-1)}}{\partial u^{(m)}\over \partial u_t},
$$
obtaining
$$
\eqalignno{{d\over dx}\biggl({\partial\hat \Lambda\over\partial
u_{tt}}\biggr)&=\sum_{j=m+1}^n \biggl[{d\over dx}\widehat{\biggl(\var {I}
{u^{(j)}(x)}\biggr)}{\partial
u^{(j-1)}\over\partial u_{t}}+\widehat{\biggl(\var {I}
{u^{(j)}(x)}\biggr)}{\partial
u^{(j)}\over\partial u_{t}}-\widehat{\biggl(\var {I}
{u^{(j)}(x)}\biggr)}{\partial
u^{(m)}\over\partial u^{(m-1)}}{\partial
u^{(m)}\over\partial u_{t}}\biggr]\cr
&=\biggl[{d\over dx}\widehat{\biggl(\var {I}
{u^{(m+1)}(x)}\biggr)}\biggr]{\partial
u^{(m)}\over\partial u_{t}}+\sum_{j=m+1}^n\biggl[
\widehat{\biggl(\var {I}
{u^{(j)}(x)}\biggr)}+{d\over dx}\widehat{\biggl(\var {I}
{u^{(j+1)}(x)}\biggr)}\biggr]{\partial
u^{(j)}\over\partial u_{t}}+\cr
&-\sum_{j=m+1}^n\biggl[\widehat{\biggl(\var {I}
{u^{(j)}(x)}\biggr)}{\partial
u^{(m)}\over\partial u^{(m-1)}}{\partial
u^{(m)}\over\partial u_{t}}\biggr]=\cr
&=\biggl[{d\over dx}\widehat{\biggl(\var {I}
{u^{(m+1)}(x)}\biggr)}\biggr]{\partial
u^{(m)}\over\partial u_{t}}+\sum_{j=m+1}^n
\biggl[\widehat{\biggl({\partial L\over\partial u^{(j)}}\biggr)}{\partial
u^{(j)}\over\partial u_{t}}-\widehat{\biggl(\var {I} {u^{(j)}(x)}\biggr)}{\partial
u^{(m)}\over\partial u^{(m-1)}}{\partial
u^{(m)}\over\partial u_{t}}\biggr].\cr}
$$
Inserting in ($a$.9) and equating it  to ($a$.8) we obtain
$$
\biggl[\widehat{\biggl({\partial L\over\partial u^{(m)}}\biggr)}-
{d\over dx}\widehat{\biggl(\var {I}
{u^{(m+1)}(x)}\biggr)}\biggr]
{\partial u^{(m)} \over \partial u_{t}}={\partial\hat
\Lambda\over\partial u^{(m-1)}_t}
{\partial
u^{(m)}\over\partial u_{t}}-\widehat{\biggl(\var {I} {u^{(j)}(x)}\biggr)}{\partial
u^{(m)}\over\partial u^{(m-1)}}{\partial
u^{(m)}\over\partial u_{t}}.
$$
But the term ${\partial
u^{(m)}\over\partial u_{t}}$ is nonzero by definition, so that
$$
{\partial\hat
\Lambda\over\partial u^{(m-1)}_t}=\widehat{\biggl(\var {I}
{u^{(m)}(x)}\biggr)}+\sum_{j=m+1}^n\widehat{\biggl(\var {I} {u^{(j)}(x)}\biggr)}
{\partial
u^{(j-1)}\over\partial u^{(m-1)}},
$$ 
which is the first step of the recurrence ($a$.7), and so on.
\bigskip
\bigskip
$\bullet$ Finally, since  ($a$.1) for $i> n-m$ coincides with ($a$.7), it
remains to prove it for $i\leq n-m$. These relations can be obtained from ($a$.2)
and ($a$.7) together with the identity  
$$
{\partial\over\partial u^{(i)}_{t}}\widehat{\biggl({d L\over dt}\biggr)}=
{\partial \over\partial u^{(i)}_{t}}\biggl({d\hat \Lambda\over dx}\biggr),
\qquad i=1,\ldots,n-m.
$$
Indeed, starting from the index $i=n-m$ and using (3.9), one may write
$$
\eqalignno{{\partial\over\partial u^{(n-m)}_{t}}\widehat{\biggl({d L\over
dt}\biggr)}&= {d\over dt}\biggl({\partial\hat L\over \partial
u^{(n-m)}_t}\biggr)+ {\partial\hat L\over \partial u^{(n-m)}}=\cr
&={d\over dt}\biggl({\partial\hat L\over \partial
u^{(n)}}
{\partial u^{(n)} \over \partial u^{(n-m)}_{t}}\biggr)+
\widehat{\biggl({\partial L\over \partial
u^{(n-m)}}\biggr)}+\sum_{j=m}^n\widehat{\biggl({\partial L\over \partial
u^{(j)}}\biggr)}
{\partial
u^{(j)}\over\partial u^{(n-m)}}.&(a.10)\cr}
$$
On the other hand
$$
{\partial\over\partial u^{(n-m)}_{t}}\biggl({d\hat \Lambda\over dx}\biggr)=
{d\over dx}\biggl({\partial\hat \Lambda\over\partial u^{(n-m)}_{t}}\biggr)+
{\partial \Lambda\over \partial u^{(n-m-1)}_{t}}+
{\partial \Lambda\over \partial u^{(m-1)}_{t}}.
$$
Performing the same steps as in the previous case, one obtains
$$
{\partial\hat \Lambda\over \partial u^{(n-m-1)}_{t}}-{d\over dt}\biggl(
{\partial\hat \Lambda\over \partial u^{(n-m-1)}_{tt}}\biggr)=
\widehat{\biggl(\var {I}
{u^{(n-m)}(x)}\biggr)}+\sum_{j=m+1}^n\widehat{\biggl(\var {I} {u^{(j)}(x)}\biggr)}
{\partial
u^{(j-1)}\over\partial u^{(n-m-1)}}
$$
which is the first recursive step of ($a$.1).

\hfill{Q.E.D.}
\bigskip
\bigskip
\bigskip
\noindent{\bf  Proof of Lemma 3.2}: We prove the Lemma by mean of the
equivalence $$
\widehat{\biggl(\var I {u(x)}\biggr)}
=\biggl[{\partial\hat \Lambda\over \partial u^{(m-1)}}-{d\over
dt}\biggl( {\partial\hat \Lambda\over \partial u^{(m-1)}_{t}}\biggr)\biggr]
{\partial
u^{(m)}\over\partial u_t}
$$
which follows from the identity
$$
{\partial\over\partial u_t}\biggl({d\hat L\over dt}\biggr)=
{\partial\over\partial u_t}\biggl({d\hat \Lambda\over dx}\biggr).
$$
Indeed, expanding, one has
$$
\eqalignno{{\partial\over\partial u_t}\biggl({d\hat L\over dt}\biggr)&=
{d\over dt}{\partial\hat  L\over\partial u_t}+
{\partial\hat L\over \partial
u}=\cr
&={d\over dt}\sum_{j=m}^n\widehat{\biggl({\partial L\over \partial
u^{(j)}}\biggr)}{\partial
u^{(j)}\over\partial u_t}+\widehat{\biggl({\partial L\over \partial
u}\biggr)}+\sum_{j=m}^n\widehat{\biggl({\partial L\over \partial
u^{(j)}}\biggr)}{\partial
u^{(j)}\over\partial u}=\cr
&={d\over dt}\biggl[{d\over dx}
\biggl({\partial \Lambda\over\partial u_{tt}}\biggr)
+{\partial \Lambda\over\partial u^{(m-1)}_{t}}{\partial u^{(m)}\over 
\partial u_t}
\biggr]+\widehat{\biggl({\partial L\over \partial
u}\biggr)}+\sum_{j=m}^n\widehat{\biggl({\partial L\over \partial
u^{(j)}}\biggr)}{\partial
u^{(j)}\over\partial u},\cr}
$$
where the last equality follows from the equivalence of ($a$.8) and ($a$.9).
Expanding the right hand side
$$
\eqalignno{{\partial\over\partial u_t}\biggl({d\hat L\over dt}\biggr)&=
{d\over dx}\biggl[{d\over dt}\biggl(
{\partial \Lambda\over\partial u_{tt}}\biggr)\biggr]+
\biggl[{d\over dt}\biggl(
{\partial \Lambda\over\partial u^{(m-1)}_{t}}\biggr)\biggr]{\partial
u^{(m)}\over  \partial u_t}+{\partial \Lambda\over\partial u^{(m-1)}_{t}}
{d\over dt}\biggl({\partial u^{(m)}\over 
\partial u_t}\biggr)+\cr
&+\widehat{\biggl({\partial L\over \partial
u}\biggr)}+\sum_{j=m}^n\widehat{\biggl({\partial L\over \partial
u^{(j)}}\biggr)}{\partial
u^{(j)}\over\partial u},&(a.11)
\cr}
$$
On the other hand 
$$
\eqalignno{{\partial\over\partial u_t}\biggl({d\Lambda\over dx}\biggr)&=
{d\over dx}\biggl({\partial \Lambda\over\partial u_t}\biggr)+
{\partial \Lambda\over \partial u^{(m-1)}}{\partial u^{(m)}\over 
\partial u_t}+{\partial \Lambda\over \partial u^{(m-1)}_t}{\partial u^{(m)}\over 
\partial u}+{\partial \Lambda\over \partial u^{(m-1)}_t}
{d\over dt}\biggl({\partial u^{(m)}\over 
\partial u_t}\biggr)=\cr
&={d\over dx}\biggl({\partial \Lambda\over\partial u_t}\biggr)+
{\partial \Lambda\over \partial u^{(m-1)}}{\partial u^{(m)}\over 
\partial u_t}+\cr
&+\biggl[\widehat{\biggl(\var I {u^{(m)}(x)}\biggr)}+
\sum_{j=m+1}^n\widehat{\biggl(\var I {u^{(j)}(x)}\biggr)}
{\dd u^{(j-1)}\over\dd u^{(m-1)}}\biggr]
{\partial u^{(m)}\over 
\partial u}+{\partial \Lambda\over
\partial u^{(m-1)}_t} {d\over dt}\biggl({\partial u^{(m)}\over 
\partial u_t}\biggr).&(a.12)\cr}
$$
Comparing ($a$.11) and ($a$.12) one obtains
$$
\eqalignno{&{d\over dx}\biggl[\biggl({\partial \Lambda\over\partial
u_t}\biggr)-
{d\over dt}\biggl(
{\partial \Lambda\over\partial u_{tt}}\biggr)\biggr]+\biggl[
\biggl({\partial \Lambda\over\partial
u^{(m-1)}}\biggr)-
{d\over dt}\biggl(
{\partial \Lambda\over\partial u^{(m-1)}_{t}}\biggr)\biggr]
{\partial u^{(m)}\over 
\partial u_t}
=\cr
&=\widehat{\biggl({\partial L\over \partial
u}\biggr)}+\sum_{j=m}^n\widehat{\biggl({\partial L\over \partial
u^{(j)}}\biggr)}{\partial
u^{(j)}\over\partial u}-
\biggl[\widehat{\biggl(\var I {u^{(m)}(x)}\biggr)}+
\sum_{j=m+1}^n\widehat{\biggl(\var I {u^{(j)}(x)}\biggr)}
{\dd u^{(j-1)}\over\dd u^{(m-1)}}\biggr]{\partial u^{(m)}\over 
\partial u}.&(a.13)\cr}
$$
But Lemma 3.1 states that
$$
\biggl[\biggl({\partial \Lambda\over\partial
u_t}\biggr)-
{d\over dt}\biggl(
{\partial \Lambda\over\partial u_{tt}}\biggr)\biggr]=
\widehat{\biggl(\var I {u_x(x)}\biggr)}-
\sum_{j=m+1}^n\widehat{\biggl(\var I {u^{(j)}(x)}\biggr)}
{\dd u^{(j-1)}\over\dd u},
$$
so that 
$$
\eqalign{{d\over dx}\biggl[
\biggl({\partial \Lambda\over\partial
u_t}\biggr)-
{d\over dt}\biggl(
{\partial \Lambda\over\partial u_{tt}}\biggr)
\biggr]&={d\over dx}\widehat{\biggl(\var I {u_x(x)}\biggr)}-{d\over dx}\biggl[
\sum_{j=m+1}^n\widehat{\biggl(\var I {u^{(j)}(x)}\biggr)}\biggr]
{\dd u^{(j-1)}\over\dd u}+\cr
&+\sum_{j=m+1}^n\widehat{\biggl(\var I {u^{(j)}(x)}\biggr)}{d\over dx}\biggl(
{\dd u^{(j-1)}\over\dd u}\biggr)=\cr
&={d\over dx}\widehat{\biggl(\var I {u_x(x)}\biggr)}
-{d\over dx}\biggl[
\sum_{j=m+1}^n\widehat{\biggl(\var I {u^{(j)}(x)}\biggr)}\biggr]
{\dd u^{(j-1)}\over\dd u}+\cr
&+\sum_{j=m+1}^n\widehat{\biggl(\var I {u^{(j)}(x)}\biggr)}\biggr(
{\dd u^{(j)}\over \dd u}-{\dd u^{(j-1)}\over \dd u^{(m-1)}}
{\dd u^{(m-1)}\over \dd u}\biggr).
\cr}
$$
Substituting in ($a$.13) we get
$$
\eqalignno{\biggl[
\biggl({\partial \Lambda\over\partial
u^{(m-1)}}\biggr)-
{d\over dt}\biggl(
{\partial \Lambda\over\partial u^{(m-1)}_{t}}\biggr)\biggr]
{\partial u^{(m)}\over 
\partial u_t}
&=\biggl[\widehat{\biggl({\partial L\over \partial
u}\biggr)}-{d\over dx}\widehat{\biggl(\var I {u_x(x)}\biggr)}\biggr]+\cr
&+\sum_{j=m}^n\biggl[\widehat{\biggl({\partial L\over \partial
u^{(j)}}\biggr)}-
{d\over dx}\widehat{\biggl(\var I {u^{(j+1)}(x)}\biggr)}\biggr]
{\partial
u^{(j)}\over\partial u}+\cr
&-\sum_{j=m+1}^n\widehat{\biggl(\var I {u^{(j)}(x)}\biggr)}\biggr(
{\dd u^{(j)}\over \dd u}-{\dd u^{(j-1)}\over \dd u^{(m-1)}}
{\dd u^{(m-1)}\over \dd u}\biggr)+\cr
&+\biggl[\widehat{\biggl(\var I {u^{(m)}(x)}\biggr)}-
\sum_{j=m+1}^n\widehat{\biggl(\var I {u^{(j)}(x)}\biggr)}
{\dd u^{(j-1)}\over\dd u^{(m-1)}}\biggr]{\partial u^{(m)}\over 
\partial u}.\cr}
$$
All the terms cancels but
$$
\biggl[
\biggl({\partial \Lambda\over\partial
u^{(m-1)}}\biggr)-
{d\over dt}\biggl(
{\partial \Lambda\over\partial u^{(m-1)}_{t}}\biggr)\biggr]
{\partial u^{(m)}\over 
\partial u_t}=\biggl[\widehat{\biggl({\partial L\over \partial
u}\biggr)}-{d\over dx}\widehat{\biggl(\var I {u_x(x)}\biggr)}\biggr]
$$

\hfill{Q.E.D.}
\bigskip
\bigskip
\noindent {\bf Proof of Lemma 3.3}: Relation (3.13) follows
from the identities
$$
{\partial\over\partial u^{(i)}}\biggl({d\hat L\over dt}\biggr)=
{\partial\over\partial u^{(i)}}\biggl({d\hat\Lambda\over dx}\biggr)
\qquad i=1,\ldots,m-1,\eqno{(a.14)}
$$
and
$$
{\partial\over\partial u^{(i)}_t}\biggl({d\hat L\over dt}\biggr)=
{\partial\over\partial u^{(i)}_t}\biggl({d\hat\Lambda\over dx}\biggr)
\qquad i=1,\ldots,m-1.\eqno{(a.15)}
$$
Starting from ($a$.14), we can write
$$
{\partial\over\partial u^{(i)}}\biggl({d\hat L\over dt}\biggr)=
{d\over dt}\biggl({\partial\hat L\over\partial u^{(i)}}\biggr)=
{d\over dt}\sum_{j=m}^n\widehat{\biggl({\partial L\over \partial
u^{(j)}}\biggr)}{\partial
u^{(j)}\over\partial u^{(i)}}+{d\over dt}
\widehat{\biggl({\partial L\over \partial
u^{(i)}}\biggr)},
$$
and
$$
{\partial\over\partial u^{(i)}}\biggl({d\hat\Lambda\over dx}\biggr)=
{d\over dx}\biggl({\partial\hat \Lambda\over\partial u^{(i)}}\biggr)+
{\partial\hat \Lambda\over \partial u^{(i-1)}}+
{\partial\hat \Lambda\over \partial u^{(m-1)}}{\partial u^{(m)}\over 
\partial u^{(i)}}+
{\partial\hat \Lambda\over \partial u^{(m-1)}_t}{d\over dt}\biggl({\partial
u^{(m)}\over  \partial u^{(i)}}\biggr).
$$
which give
$$
{d\over dt}\biggl[\sum_{j=m}^n\widehat{\biggl({\partial L\over \partial
u^{(j)}}\biggr)}{\partial
u^{(j)}\over\partial u^{(i)}}+
\widehat{\biggl({\partial L\over \partial
u^{(i)}}\biggr)}\biggr]={d\over dx}\biggl({\partial\hat \Lambda\over\partial
u^{(i)}}\biggr)+ {\partial\hat \Lambda\over \partial u^{(i-1)}}+
{\partial\hat \Lambda\over \partial u^{(m-1)}}{\partial u^{(m)}\over 
\partial u^{(i)}}+
{\partial\hat \Lambda\over \partial u^{(m-1)}_t}{d\over dt}\biggl({\partial
u^{(m)}\over  \partial u^{(i)}}\biggr).\eqno{(a.16)}
$$
\bigskip
\bigskip
\noindent On the other hand, in  ($a$.15)
$$
{\partial\over\partial u^{(i)}_t}\biggl({d\hat L\over dt}\biggr)
{\partial\hat L\over\partial u^{(i)}}+{d\over dt}\biggl(
{\partial\hat L\over\partial u^{(i)}_t}\biggr)
=\widehat{\biggl({\partial L\over \partial
u^{(i)}}\biggr)}+
\sum_{j=m}^n\widehat{\biggl({\partial L\over \partial
u^{(j)}}\biggr)}{\partial
u^{(j)}\over\partial u^{(i)}}+{d\over dt}
\biggl({\partial\hat L\over \partial
u^{(i)}_t}\biggr)
$$
and 
$$
{\partial\over\partial u^{(i)}_t}\biggl({d\hat\Lambda\over dx}\biggr)=
{d\over dx}\biggl({\dd\Lambda\over \dd u^{(i)}_t}\biggr)+
{\dd\Lambda\over \dd u^{(i-1)}_t}+{\dd\Lambda\over \dd u^{(m-1)}_t}
{\partial
u^{(m)}\over\partial u^{(i)}}
$$
which give
$$
\widehat{\biggl({\partial L\over \partial
u^{(i)}}\biggr)}+
\sum_{j=m}^n\widehat{\biggl({\partial L\over \partial
u^{(j)}}\biggr)}{\partial
u^{(j)}\over\partial u^{(i)}}=-{d\over dt}\biggl(
{\partial\hat L\over\partial u^{(i)}_t}\biggr)+{d\over
dt}\biggl({\dd\Lambda\over \dd u^{(i)}_t}\biggr)+ {\dd\Lambda\over \dd
u^{(i-1)}_t}+{\dd\Lambda\over \dd u^{(m-1)}_t} {\partial
u^{(m)}\over\partial u^{(i)}}.
$$
Performing the derivative w.r.t. $t$ and substituting into ($a$.16) one gets
$$
\eqalign{&-{d^2\over dt^2}\biggl(
{\partial\hat L\over\partial u^{(i)}_t}\biggr)+{d\over dx}{d\over
dt}\biggl({\dd\Lambda\over \dd u^{(i)}_t}\biggr)+ {d\over dt}{\dd\Lambda\over
\dd u^{(i-1)}_t}+{d\over dt}\biggl[{\dd\Lambda\over \dd u^{(m-1)}_t} {\partial
u^{(m)}\over\partial u^{(i)}}\biggr]=\cr
&={d\over dx}\biggl({\partial\hat \Lambda\over\partial
u^{(i)}}\biggr)+ {\partial\hat \Lambda\over \partial u^{(i-1)}}+
{\partial\hat \Lambda\over \partial u^{(m-1)}}{\partial u^{(m)}\over 
\partial u^{(i)}}+
{\partial\hat \Lambda\over \partial u^{(m-1)}_t}{d\over dt}\biggl({\partial
u^{(m)}\over  \partial u^{(i)}}\biggr),\cr}
$$
which gives
$$
\eqalignno{&-{d^2\over dt^2}\biggl(
{\partial\hat L\over\partial u^{(i)}_t}\biggr)={d\over dx}\biggl[
\biggl({\partial\hat \Lambda\over\partial
u^{(i)}}\biggr)-
{d\over
dt}\biggl({\dd\Lambda\over \dd u^{(i)}_t}\biggr)\biggr]+
\biggl[{\partial\hat \Lambda\over \partial u^{(i-1)}}
-{d\over dt}{\dd\Lambda\over
\dd u^{(i-1)}_t}\biggr]+\cr
&+\biggl[
{\partial\hat \Lambda\over \partial u^{(m-1)}}
-{d\over dt}{\dd\Lambda\over \dd u^{(m-1)}_t}\biggr] {\partial
u^{(m)}\over\partial u^{(i)}}.&(a.17)
\cr}
$$
The left hand side of ($a$.17) is zero if $i>n-m$.

\noindent If $i<n-m$,
$$
{\partial\hat L\over \partial
u^{(i)}_t}= \sum_{k=m+i}^n \widehat{\biggl(
{\partial L\over \partial u^{(k)}}\biggr)}{\dd u^{(k)}\over\dd
u^{(i)}_t}=\biggr[
{\dd\hat\Lambda \over\dd u^{(i-1)}_{tt}}+{d\over dx}
\biggl({\dd\hat\Lambda \over\dd u^{(i)}_{tt}}\biggr)\biggr].\eqno{(a.18)}
$$
Indeed, using ($a$.1), and relation (3.6b), which we rewrite here below,
$$
{d\over dx}\biggl({\dd u^{(k)}\over\dd
u^{(i)}_t}\biggr)={\dd u^{(k+1)}\over\dd
u^{(i)}_t}-{\dd u^{(k)}\over\dd
u^{(i-1)}}
$$
one obtains
$$
\eqalign{{\dd\hat\Lambda \over\dd u^{(i-1)}_{tt}}&+{d\over dx}
\biggl({\dd\hat\Lambda \over\dd u^{(i)}_{tt}}\biggr)=
\sum_{k=i+m}^n \widehat{\biggl(
\var I {u^{(k)}(x)}\biggr)}{\dd u^{(k-1)}\over\dd
u^{(i-1)}_t}+{d\over dx}\biggl[\sum_{k=i+m+1}^n \widehat{\biggl(
\var I {u^{(k)}(x)}\biggr)}{\dd u^{(k-1)}\over\dd
u^{(i)}_t}\biggr]=\cr
&=\sum_{k=i+m+1}^n\biggl[{d\over dx}\widehat{\biggl(
\var I {u^{(k)}(x)}\biggr)}\biggr]{\dd u^{(k-1)}\over\dd
u^{(i)}_t}+\sum_{k=i+m+1}^n\widehat{\biggl(
\var I {u^{(k)}(x)}\biggr)}{\dd u^{(k)}\over\dd
u^{(i)}_t}+\widehat{\biggl(
\var I {u^{(i+m)}(x)}\biggr)}{\dd u^{(i+m-1)}\over\dd
u^{(i-1)}_t}=\cr
&=\sum_{k=i+m+1}^n\widehat{\biggl(
{\dd L\over u^{(k)}}\biggr)}{\dd u^{(k)}\over\dd
u^{(i)}_t}+\biggl[{d\over dx}\widehat{\biggl(
\var I {u^{(m+i+1)}(x)}\biggr)}\biggr]{\dd u^{(m+i)}\over\dd
u^{(i)}_t}=\cr
&= \sum_{k=m+i}^n \widehat{\biggl(
{\partial L\over \partial u^{(k)}}\biggr)}{\dd u^{(k)}\over\dd
u^{(i)}_t}
\cr}
$$
where the last identity follows from  the fact that
$$
{\dd u^{(k)}\over\dd
u^{(l)}_t}=0
$$
if $k-l<m$. In (3.6), this implies
$$
{\dd u^{(k+m)}\over\dd
u^{(k)}_t}={\dd u^{(k+m+j)}\over\dd
u^{(k+j)}_t}.
$$
Finally, substituting ($a$.18) into ($a$.17) gives (3.13)  
$$
\eqalignno{\biggl[
{\partial\hat \Lambda\over \partial u^{(i-1)}}&-{d\over dt}\biggl(
{\partial\hat \Lambda\over \partial u^{(i-1)}_{t}}\biggr)
+{d^2\over dt^2}\biggl(
{\partial\hat \Lambda\over \partial u^{(i-1)}_{tt}}\biggr)\biggr]=\cr
&=\biggl[{\partial\hat \Lambda\over \partial u^{(m-1)}}-{d\over dt}\biggl(
{\partial\hat \Lambda\over \partial u^{(m-1)}_{t}}\biggr)\biggr]
{\partial u^{(m)}\over  \partial u^{(i)}}-{d\over dx}\biggl[
{\partial\hat \Lambda\over \partial u^{(i)}}-{d\over dt}\biggl(
{\partial\hat \Lambda\over \partial u^{(i)}_{t}}\biggr)
+{d^2\over dt^2}\biggl(
{\partial\hat \Lambda\over \partial u^{(i)}_{tt}}\biggr)\biggr].\cr}
$$
\hfill{Q.E.D.}

\bigskip
\bigskip
\noindent {\scap 4.  Applications to Painlev\'e equations}
\bigskip
\bigskip
In this Section we  study some applications of Theorem 2.1 and Theorem 3.1.
to  show how
the finite dimensional Hamiltonian structure of Painlev\' e equations comes from
an infinite dimensional structure via the above procedure.
\bigskip
\noindent {\scap 4.1 PI as scaling reduction of KdV}
\bigskip
At the beginning we study the problem following the Hamiltonian scheme,
then we will apply the framework of Theorem 3.1.

\bigskip
We consider the KdV equation 
$$
u_t=6uu_x-u_{xxx}.\eqno{(4.1)}
$$
i.e. the $t=t_{_{1}}$--flow in the KdV hierarchy (1.2); it  
admits the  nondegenerate {\it scaling} symmetry
$$
I_{_{(s)}}=\int (u^3+{u_x^2\over 2}+2u x+6t u^2)dx,\eqno{(4.2)}
$$
which depends on $x,u,u_x,t$.
We note that $L=[\pd 1 L+4x \pd {-1} L+12t \pd 0 L]$,where $\pd {-1} I=\int
\pd {-1} L dx=\int {u(x)\over 2}dx$ and $\pd {0} I=\int \pd {0} L dx=\int
{u^2(x)\over 2}dx$ are the first Hamiltonians of the KdV hierarchy.
\bigskip
\noindent
Theorem 2.1 states that the $t$--flow is
Hamiltonian on the manifold $\S$ of the stationary points of the symmetry, i.e.
$\S$  is the $2$--dimensional manifold of the
solutions of the Euler--Lagrange equation
$$
\var {I}  {u(x)}=u_{xx}-3u^2-2x-12t u=0.\eqno{(4.3)}
$$

\noindent It is invariant under the $t$--th flow and it naturally carries the
system of canonical coordinates:
 $$
\cases{q=u\cr 
p=\var  I  {u_x}=u_x\cr}
$$ 
Notice that  the identities
$$
\cases{p_x+{\dd  H\over{\dd q}}\equiv -
\var { I}  {u}\cr
q_x-{\dd { H}\over{\dd p}}\equiv 0,\cr}\eqno{(4.3)}
$$
hold, where $H$ is the generalized Legendre transform of $L$:
$$
H=-L+ \var  I  {u_x}
u_x.
$$
The first of  identities (4.4) allows us to express the higher derivatives
$u^{(m)}$ for $m\geq 2$ in terms of $x$,$t$, $p,\ q$ and
$p^{{_{(l)}}}$ with $l=1,\dots,{m-2+1}$.
\bigskip
\noindent
On $\S$  $p_x+{\dd { H}\over{\dd q_{_1}}}\equiv
0$, and the system (4.4) reduces to the canonical Hamiltonian system
$$
\cases{p_x=3q^2+2x+12t q\cr
q_x=p,\cr}
$$
for the 
Hamiltonian  function
$$
 H=-L+u_x^2={p^2\over 2}-q^3-2q x-6t q^2\eqno{(4.5)}
$$
 giving the reduced $x$--flow. This system is equivalent to the second order
ODE in the variable $q$ :
$$
q^{''}=-3q^2-2x-12t q.\eqno{(4.6)}
$$
\bigskip
\noindent The space $\S$ is the set of the stationary points of the scaling
symmetry (4.2); this means that $\S$ carries a ``natural" system of canonical
coordinates $\{w_i,\pi_i\}$, given by the self--similar function of $u$, i.e;
combinations of $u,x,t$ in the variable $z (x,t)$ invariant w.r.t. the scaling.
We will call them scaling coordinates. In this case
$$
\cases{w={q\over 2}+t\cr
\pi=2p\cr}
$$
with $z=x-6t^2$.

In terms of the scaling coordinates
the system reads
  $$
\cases{{d\pi\over dz}=-{\dd \H\over\dd w}\cr
{dw\over dz}={\dd \H\over\dd \pi}.\cr}
$$
for the Hamiltonian
$$
\H={\pi^2\over 8}-8w^3-4wz+8t^3+4tz.
$$
The system is equivalent to the ODE:
$$
w^{''}=6w^2+z,\eqno{(4.7)}
$$
that is exactly Painlev\'e I. The Hamiltonian $\H$ differs from the usual PI
Hamiltonian for the terms in $z,t$  that do not enter in the
Hamiltonian system. 
\bigskip 
\noindent We now construct the 
time dependent Hamiltonian function $(-\tilde Q)$, that is the
reduction on $\S$ of 
$$
-Q=-\Lambda +p {d q\over dt},
$$
 where $p,q$ are  expressed in terms
of $u,u_x$,  and $\Lambda(x,t,u,u_x,u_{xx},u_{xxx})$, calculated from 
$$
{dL\over dt}={d\Lambda\over dx}.
$$
has the form
$$
\Lambda=6t(4u^3-2uu_{xx}+u_x^2)+2x(3u^2-u_{xx})+
{9\over 2}u^4+{1\over 2}u_{xx}^2+2u_x-3u^2u_{xx}+6uu_x^2-u_xu_{xxx}.
$$

By direct calculation one obtains

$$
Q=12t(2u^3+{u_x^2\over 2}-uu_{xx})+{u_{xx}^2\over
2}-3u^2u_{xx}+{9\over 2}u^4+2u_x+2x(3u^2-u_{xx}), \eqno{(4.8)}
$$
This reduces on $\S$ to
$$
\tilde Q=12t({p^2\over 2}-q^3-6tq^2-2xq)+2p-2x^2
\eqno{(4.9)}
$$
Theorem 2.1 states that $(-\tilde Q)$ is the Hamiltonian for the
reduced $t$-flow,i.e., in terms of $p$ and $q$
 $$
\cases{\dot q=-2\ (6t p+1)=-{\dd\tilde Q\over\dd p}\cr
\dot p=-12t\ (3q^2+2x+12tq)={\dd\tilde Q\over\dd
q},\cr}\eqno{(4.10)} $$
Notice that  system (4.10), written in terms of the scaling coordinates $w$
and $z$, gives the same Painlev\' e I.
\bigskip
\noindent{\bf Remark}: In this case the evolution equation is
Hamiltonian and  it can be written in the form
$$
u_t=\pois{u(x),I_1}={d\over dx} \var {\pd 1 I}  u,
$$
where
$
\pd 1 I=\int \pd 1 L dx
\ $ with  density
$$
\pd 1 L=u^3+{u_x^2\over 2}.
$$
On the other hand the scaling symmetry defines the stationary flow 
$$
{d u\over ds}=12tu_x+u_t+2=0.
$$
which is  Hamiltonian:
$$
{d u\over ds}=\pois{u(x),I}={d\over dx} \var  I 
u=0.
$$

The $s$-flow and the $t$-flow commute, but the Hamiltonian generating the
scaling depends explicitly on the time $t$, so that the relation
$$
{d\pd {(s)} I\over dt}=\pois{\pd {(s)} I,\pd 1 I}+{\dd \pd {(s)} I\over \dd t}=0
$$
holds. Hence we have an alternative way to define the reduced Hamiltonian
$Q$, following [BN]:
$$
 {d
\over
dx}Q
={d L\over d t}-{\partial I\over\partial u}{d\over dx}
{\partial I_1\over\partial u}.
$$
In this case the relation (4.2) follows as a consequence.
\bigskip
\bigskip
\bigskip
\noindent System (4.10) i.e. the reduction of the
$t$--flow on $\S$, can be obtained from the Lagrangian point of view; indeed
one can consider the evolution equation (4.1) as the definition of $u_{xxx}$ in
terms of $(u, u_x,u_{xx},u_t)$, explicitly: 
 $$
 u_{xxx}=6uu_x-u_t. 
$$
Differentiating this relation in $x$ one obtains 
$$
\cases{u^{(4)}=6uu_{xx}+6u_x^2-u_{xt}\cr
u^{(5)}=18u_xu_{xx}+36u^2-6uu_t-u_{xxt}\cr
\vdots\cr}
$$
which is a map from the ''coordinates"
$$
u, u_x,u_{xx}, u^{(3)},u^{(4)},u^{(5)},u^{(6)},\ldots
$$
into
$$
u, u_x,u_{xx},u_t, u_{xt},u_{xxt},u_{tt},\ldots
$$
\bigskip 
\noindent  The Lagrangian $L$ depends on $
u, u_x$ and hence its derivative ${dL\over dt}$ looks like
$$
{dL\over dt}={\partial L\over \partial t}+{\partial L\over \partial u}
u_t+{\partial L\over \partial u_x}
u_{xt}.
$$
And there exist a functional 
$\hat\Lambda$ depending on $u, u_x,u_{xx},u_t$ such that
$$
{d L\over dt}={d\hat\Lambda\over dx}={\partial\hat \Lambda\over \partial
x}+{\partial\hat \Lambda\over \partial u} u_x
+{\partial\hat \Lambda\over \partial u_x} u_{xx}+{\partial\hat \Lambda\over
\partial u_{xx}} u^{(3)}+{\partial\hat
\Lambda\over \partial u_t} u_{xt}.
$$
Here,
in terms of the new coordinates
$$
\hat\Lambda=6t(4u^3-2uu_{xx}+u_x^2)+2x(3u^2-u_{xx})+
{9\over 2}u^4+{1\over 2}u_{xx}^2+2u_x-3u^2u_{xx}+u_xu_t,
$$

\noindent$\hat L\equiv L$, and
$$
\eqalign{{\partial L\over \partial t}&=6u^2\cr
{\partial L\over \partial u}&=3u^2+2x+12tu\cr
{\partial L\over \partial u_x}&=u_x\cr}
$$
\bigskip
$$
\eqalign{{\partial\hat \Lambda\over \partial x}&=6u^2-2u_{xx}\cr
{\partial\hat\Lambda\over \partial u}&=72tu^2-12tu_{xx}+12xu+18u^3-6uu_{xx}\cr
{\partial\hat \Lambda\over \partial u_x}&=12tu_x+u_t+2\cr
{\partial\hat \Lambda\over \partial u_{xx}}&=-12tu+u_{xx}-2x-3u^2\cr
{\partial\hat \Lambda\over \partial u_t}&=u_x\cr}
$$
\bigskip
\noindent The condition $\var I u=0$, that
defines the submanifold $\S$, is equivalent to the condition
$$
{\partial\hat \Lambda\over \partial u_{xx}}=0
$$ 
Hence we have an alternative definition of the space $\S$, and an alternative
way to defines the canonical coordinates:
$$
\cases{q=u\cr
p={\partial\hat \Lambda\over \partial u_t}=u_x\cr}
$$
\bigskip
\noindent Theorem 3.1 states that the reduced $t$--flow is Lagrangian ,with
Lagrangian $\Lambda$, in this case it is easy to verify it, indeed, on $\S$,

$$
\cases{{\partial \hat\Lambda\over \partial u_{xx}}=-12tu+u_{xx}-2x-3u^2=0\cr
{\partial \hat\Lambda\over \partial u_x}=12tu_x+u_t+2=0\cr
{\partial \hat\Lambda\over \partial u}-
{d\over dt}\bigl({\partial \Lambda\over \partial u_t}\bigr)=
72tu^2-12tu_{xx}+12xu+18u^3-6uu_{xx}-u_{xt}=0,
\cr}
$$
where the first equation is the definition of  the submanifold $\S$ itself, the
other two reproduces (4.9), indeed they can be rewritten as 
$$
\cases{u_t=-12tu_x-2\cr
u_{xt}=-12t(3u^2+2x+12tu) \cr}
$$ 
\bigskip
\bigskip
\noindent {\scap 4.2  PII as scaling reduction of mKdV}
\bigskip
\noindent One can repeat the same procedure as in section 4.1 starting from the 
mKdV equation  
$$
u_t=6u^2u_x-u_{xxx}.\eqno{(4.11)}
$$ 
It  
admits the  nondegenerate {\it scaling} symmetry
$$
I=\int ({3\over 2}t(u^4+u_x^2)+{u^2x\over 2})dx,\eqno{(4.12)}
$$
which depends on $x,u,u_x,t$.
We notice that $L=3t\pd 1 L+{u^2x\over 2}$.

Here
$\S$  is the $2$--dimensional manifold of the
solutions of the Euler--Lagrange equation
$$
\var {I}  {u(x)}=u_{xx}-{1\over 3t}(6tu^3+ux)=0.\eqno{(4.13)}
$$
It naturally carries the
system of canonical coordinates:
 $$
\cases{q=u\cr 
p=\var  I  {u_x}=3tu_x.\cr}
$$ 
As in the previous case we read
the Euler--Lagrange equation as a reduced $x$--flow 
 with Hamiltonian
$$
 H={p^2\over 6t}-{3\over 2}tq^4-{1\over 2}q^2 x\eqno{(4.14)}
$$ 

where $H$ is the generalized Legendre transform of $L$:
$$
H=-L+ \var  I  {u_x}
u_x=-L+3tu_x^2.
$$
The system is equivalent to the second order
ODE in the variable $q$ :

$$
q_{xx}=2q^3+{1\over 3t}qx.
$$
The scaling coordinates are now
$$
\cases{ w=
(3t)^{^{1\over 3}}q\cr
 \pi={p\over (3t)^{^{1\over 3}}}\cr}
$$
in the variable $z={x\over (3t)^{^{1\over
3}}}$,
and the system transforms into
  $$
\cases{{d\pi\over dz}=-{\dd \H\over\dd w}\cr
{dw\over dz}={\dd \H\over\dd \pi}.\cr}
$$
for the Hamiltonian
$$
\H={1\over 2}(\pi^2-w^4)-{1\over 2}zw^2.
$$
The system is equivalent to the ODE:
$$
w^{''}=2w^3+zw.\eqno{(4.15)}
$$

that is exactly Painlev\'e II.

\bigskip
\noindent We now construct the 
time dependent Hamiltonian function $(-\tilde Q)$, that is the
reduction on $\S$ of 
$$
-Q=-\Lambda +p {d q\over dt},
$$
where
$$
\Lambda=t\bigl(6u^6-6u^3u_{xx}+18u^2u_x^2+{3\over
2}u_{xx}^2-3u_xu_{xxx}\bigr)+x\bigl({3\over
2}u^4+{1\over
2}u_x^2-uu_{xx}\bigr)+ uu_x. 
$$
By direct calculation one obtains

$$
Q=6t(u^6-u^3u_{xx}+{1\over 4}u_{xx}^2)+x({3\over
2}u^4-uu_{xx}+{1\over 2}u_x^2)+uu_x,
\eqno{(4.16)} 
$$
which on $\S$ reduces to
$$
\tilde Q={x\over 2}({p^2\over 9t^2}-q^4)+{1\over
3t}pq-{1\over
6t}q^2x^2
\eqno{(4.17)}
$$
and  is the Hamiltonian for the reduced $t$-flow .
In fact
 $$
\cases{\dot q=-{q\over 3t}-{xp\over 9t^2}=-{\dd\tilde Q\over\dd
p}\cr \dot p={p\over 3t}-{qx^2\over 3t}-2q^3x= 
{\dd\tilde Q\over\dd q},\cr}\eqno{(4.18)}
$$
Notice that also the system (4.18), written in $w$ and $z$, gives Painlev\' e II.
\bigskip

\noindent{\bf Remark}: The evolution equation is
Hamiltonian and  can be written in the form
$$
u_t=\pois{u(x),I_1}={d\over dx} \var {\pd 1 I}  u,
$$
where
$
\pd 1 I=\int \pd 1 L dx
\ $ with  density
$$
\pd 1 L={1\over 2}(u^4+u_x^2)
$$
On the other hand the scaling symmetry defines the Hamiltonian stationary flow 
$$
{d u\over ds}=xu_x+3tu_t+u={d\over dx} \var {\pd {(s)} I}  u=0.
$$
\bigskip
\bigskip
\bigskip
\noindent We now deduce system (4.18)  from the
Lagrangian point of view, reading the evolution equation (4.11)
as the definition of $u_{xxx}$ in terms of $(u, u_x,u_{xx},u_t)$, explicitly: 
$$
u_{xxx}=6u^2u_x-u_t.
$$
Differentiating this relation in $x$ one obtains 
$$
\cases{u^{(4)}=6u^2u_{xx}+12uu_x^2-u_{xt}\cr
u^{(5)}=36uu_xu_{xx}+36u^4u_x+12u_x^3-6u^2u_t-u_{xxt}\cr}
$$
which is a map from the ''coordinates"
$$
u, u_x,u_{xx}, u^{(3)},u^{(4)},u^{(5)},u^{(6)},\ldots
$$
into
$$
u, u_x,u_{xx},u_t, u_{xt},u_{xxt},u_{tt},\ldots
$$
\bigskip 
\noindent Here
$$
\hat \Lambda=t\bigl(6u^6-6u^3u_{xx}+{3\over
2}u_{xx}^2-3u_xu_{xxx}+3u_xu_t\bigr)+x\bigl({3\over
2}u^4+{1\over
2}u_x^2-uu_{xx}\bigr)+ uu_x. 
$$
and
$$
\eqalign{{\partial L\over \partial t}&={3\over 2}(u^4+u_x^2)\cr
{\partial L\over \partial u}&=6tu^3+ux\cr
{\partial L\over \partial u_x}&=3tu_x\cr}
$$
\bigskip
$$
\eqalign{{\partial \Lambda\over \partial x}&={3\over 2}(u^4+u_x^2)-uu_{xx}\cr
{\partial\Lambda\over \partial u}&=36tu^5-18tu^2u_{xx}+6xu^3-xu_{xx}+u_x
\cr
{\partial \Lambda\over \partial u_x}&=xu_x+u+3tu_t\cr
{\partial \Lambda\over \partial u_{xx}}&=3tu_{xx}-6tu^3-ux\cr
{\partial \Lambda\over \partial u_t}&=3tu_{x}\cr}
$$

\bigskip
\bigskip
\noindent The condition $\var I u=0$, that
defines the submanifold $\S$, is equivalent to the condition
$$
{\partial \Lambda\over \partial u_{xx}}=3tu_{xx}-6tu^3-ux=0
$$ 
Hence we have an alternative definition of the space $\S$, and an alternative
way to defines the canonical coordinates:
$$
\cases{q=u\cr
p={\partial \Lambda\over \partial u_t}=3tu_x\cr}
$$
\bigskip
\noindent On $\S$
$$
\cases{{\partial \Lambda\over \partial u_{xx}}=3tu_{xx}-6tu^3-ux=0\cr
{\partial \Lambda\over \partial u_x}=xu_x+3tu_t+u=0\cr
{\partial \Lambda\over \partial u}-
{d\over dt}\bigl({\partial \Lambda\over \partial u_t}\bigr)=
-18tu^2u_{xx}+36tu^5+6xu^3-xu_{xx}-2u_x-3tu_{xt}=0,
\cr}
$$
where the first defines the submanifold, the second one gives the motion of $u$ and the third 
the motion of $u_x$, hence one can rewrite them as 
$$
\cases{3tu_t=-xu_x-u\cr
3tu_{xt}=u_x-2u^3x-{x^2u\over 3t},\cr}
$$
which coincides with (4.18). 
\bigskip
\bigskip
{\sscap 4.3  PIII as scaling reduction of Sine-Gordon}
\bigskip
A particular case of Painlev\' e III equation can be obtained as reduction of the Sine--Gordon
equation 
$$
\cases{u_t=v= \var {\pd 1 I}  v\cr
v_t=u_{xx}-\sin\  u=-\var {\pd 1 I}  u,\cr}\eqno{(4.19)}
$$
via the scaling
$$
\cases{{d u\over ds}=xv +tu_x=\var {\pd {(s)} I}  v=0\cr
{d v\over ds}=x(u_{xx}-\sin\ u)+tv_x+u_x=-\var {\pd {(s)} I} 
u=0,\cr}\eqno{(4.20)} $$

where
$
\pd 1 I=\int \pd 1 L dx
$, $\pd {(s)} I=\int L dx$, with the Hamiltonians
$$
\pd 1 L={1\over 2}(v^2+u_x^2)-\cos\ u
$$
and
$$
L={x\over 2}(v^2+u_x^2)-x \cos\ u+tvu_x=x\pd 1 L+tvu_x
$$
w.r.t. the Poisson bracket 
$$
\pois{F,G}=\int (\var
f {u(x)} \var g  {v(x)}-\var f  {v(x)} \var g  {u(x)})\ dx. 
 $$

The scaling reduction equation means
$$
\var {\pd {(s)} I}  {u(x)}=\var {\pd {(s)} I}  {v(x)}=0,
$$
which defines the submanifold $\S$:
$$
u_x=-{xv\over t},\eqno{(4.21)}
$$
with the canonical coordinates 
$$
\cases{p=xu_x+tv={t^2-x^2\over t}v\cr
q=u\cr}
$$
 in
$\S$.

The equation defining $\S$ can be written as an Hamiltonian system in
canonical form  describing  the
reduced $x$-flow:
$$
\cases{(p)_x=-{\dd H\over{\dd q}}\cr
(q)_x={\dd H\over{\dd p}}\cr}
$$
where
$$
 H={x\over 2}{p^2\over x^2-t^2}+x\ \cos\ q.\eqno{(4.22)}
$$
In terms of the scaling coordinates
$$
\cases{w=q\cr
 \pi={-2z\over t}p\cr}
$$
in the  variable $z={x^2-t^2\over
2}$, 
the Hamiltonian system transform into
 $$
\cases{{d\pi\over dz}=-{\dd \H\over\dd w}\cr
{dw\over dz}={\dd \H\over\dd \pi}.\cr}
$$
for the Hamiltonian
$$
\H=-{1\over 4}{\pi^2\over z}-\cos\ w.
$$
The system is equivalent to Painlev\'e III:
$$
2zw^{''}+2w^{'}-\sin\ w=0.
$$  
\bigskip
\noindent Let us now construct the 
time dependent Hamiltonian function $(-\tilde Q)$, that is the
reduction on $\S$ of 
$$
-Q=-\Lambda +p {d q\over dt},
$$
where
$$
\Lambda=xu_xv+t\bigl({1\over 2}(v^2+u_x^2)-\cos u\bigr). 
$$
By direct calculation one obtains
$$
Q=t\bigl({1\over 2}(v^2-u_x^2)-\cos u\bigr)
$$
which on $\S$ reduces to

$$
\tilde Q=-t\ ({1\over 2}{p^2\over t^2-x^2}- \cos\ q).
$$
This is the Hamiltonian for the reduced $t$-flow .
In fact
 $$
\cases{\dot q=v=-{t\over x}u_x=-t{p\over t^2-x^2}=-{\dd\tilde Q\over\dd p}\cr
\dot p=v+xv_x+tv_t=t\ \sin\ q={\dd\tilde Q\over\dd
q}.\cr} \eqno{(4.23)}
$$
Note that also the system (4.23), written in $w$ and $z$, gives Painlev\' e III.

\noindent{\bf Remark }: We now deduce system (4.23)  from the
Lagrangian point of view, reading the evolution equation (4.19)
as the definition of $v$ in terms of $u_t$, explicitly: 
$$
\cases{v=u_t\cr
u_{xx}=v_t+\sin u=u_{tt}+\sin u\cr}
$$
Differentiating this relation in $x$ one obtains 
$$
\eqalign{v_x&=u_{xt}\cr
v_{xx}&=v_{tt}-v\cos u=u_{ttt}-u_t\cos u\cr
\vdots&\cr 
u_{xxx}&=u_{ttt}+u_x cos u\cr
\vdots&\cr}
$$
which is a map from the ''coordinates"
$$
u,v, u_x,v_x, u_{xx},v_{xx},\ldots
$$
into
$$
u, u_x, u_t, u_{xt}, u_{tt},u_{xtt},\ldots
$$
\bigskip 
\noindent Here
$$
\hat \Lambda=xu_xu_t+t\bigl({1\over 2}(u_t^2+u_x^2)-\cos u\bigr) 
$$
and
$$
\hat L=tu_xu_t+x\bigl({1\over 2}(u_t^2+u_x^2)-\cos u\bigr)
$$
which give
$$
\eqalign{{\partial\hat L\over \partial t}&={1\over 2}(u_t^2+u_x^2)-\cos u\cr
{\partial\hat L\over \partial u}&=-\sin u\cr
{\partial\hat L\over \partial u_x}&=tu_t+xu_x\cr
{\partial\hat L\over \partial u_t}&=tu_x+xu_t\cr}
$$
and
$$
\eqalign{{\partial\hat \Lambda\over \partial x}&={1\over 2}(u_t^2+u_x^2)-\cos
u\cr {\partial\hat \Lambda\over \partial u}&=-t\sin u\cr
{\partial\hat \Lambda\over \partial u_x}&=xu_t+tu_x\cr
{\partial\hat \Lambda\over \partial u_t}&=xu_x+tu_t\cr}
$$

\bigskip
\bigskip
\noindent The condition $\var I u=0$, that
defines the submanifold $\S$, is equivalent to the condition
$$
{\partial\hat \Lambda\over \partial u_{x}}=xu_t+tu_x=0
$$ 
Hence we have an alternative definition of the space $\S$, and an alternative
way to defines the canonical coordinates:
$$
\cases{q=u\cr
p={\partial\hat \Lambda\over \partial u_t}=xu_x+tu_t\cr}
$$
\bigskip
\noindent On $\S$
$$
\cases{{\partial\hat \Lambda\over \partial u_{x}}=xu_t+tu_x=0\cr
{\partial\hat \Lambda\over \partial u}-
{d\over dt}\bigl({\partial\hat \Lambda\over \partial u_t}\bigr)=
-t\ \sin u-xu_{xt}-u_t-tu_{tt}=0,
\cr}
$$
where the first defines the submanifold, and the second
$$
-t\sin u-xu_{xt}-u_t-tu_{tt}=0
$$
 coincides with (4.23).

\bigskip
\bigskip
\noindent {\scap 5.  Self--similar solutions of n--waves equation and 
Hamiltonian MPDEs}
\bigskip
\bigskip
\noindent {\scap 5.1 $n$--waves equations and their symmetries}
\bigskip

\noindent Let us consider the equation
$$
u_t-v_x-[u,v]=0,\eqno{(5.1)}
$$
where
$$
u=[\l,a]\ \ \ v=[\l,b]\ \ \ a= diag\ (a^1,....a^n)
\ \ \ b= diag\ (b^1,....b^n)\eqno{(5.2)}
$$
and $\l$ is a function of $x,\ \ t$.

Following [DS] it is possible to rewrite  (5.1) as an infinite
dimensional Hamiltonian system on the space $\M$ of functions of $x$ with values
in $Mat(n,C)$ with the inner product
$$
(u,v)=\int \tr\ (u(x) v(x))dx.
$$
On the space $\F$ of functionals

$$
F=\int f(x,u,u_x,....u^{^{(k)}})\ dx
$$
one can define $\nabla_{_{u}}F \in \M$  by
$$
{d\over d\epsilon}F(u+\epsilon w)\vert_{_{\epsilon=0}}=(\nabla_{_{u}}F,w)
$$
and the Poisson structure $P$ with the Poisson bracket
$$
\pois 
{F,G}(u)=(\nabla_{_{u}}F,[\nabla_{_{u}}G,{d\over dx}+u])\eqno{(5.3)}
$$

The $n$-waves equation (5.1) is a Hamiltonian system w.r.t. this  Poisson
structure: 
$$
u_t=P dI_1
=[\nabla_{_{u}}I_1,{d\over dx}+u]=[-v,{d\over dx}+u],\eqno(5.4)
$$
where
$$
I_1=\int L_1 dx=-{1\over 2}\int Tr\ (uv)dx
$$
that in components of $\l$ gives
$$
I_1=\int \sum_i\sum_k ({b_i\over
a_k-a_i}u_{ik}u_{ki}) dx=\int \sum_i\sum_k [(b_i-b_k)
(a_i-a_k)\ga i k\ga k i] dx.\eqno{(5.5)}
$$

\bigskip

For ($n=3$, $u^{^{T}}=-u$) one can reduce to a particular
case of P VI equation (see [D], where the Hamiltonian structure for this particular case of P VI is derived from the Hamiltonian structure of the $n$--waves
equation.), imposing
the scaling  $$
{d u\over ds}=tu_t+xu_x+u=0.\eqno(5.6)
$$
It admits the Hamiltonian form
$$
{d u\over ds}=[\nabla_{_{u}}I_{(s)},{d\over dx}+u]=0\eqno(5.7)
$$
where
$$
I_{(s)}=\int L dx=-{1\over 2}\int Tr\ (tuv+xu^2)dx=\int\sum_i\sum_k (t{b_i\over
a_k-a_i}-{x\over 2})u_{ik}u_{ki}dx
$$
and
$$
\nabla_{_{u}} I_{(s)}=-tv-xu,\eqno(5.8)
$$
or, in terms of $\l$:
$$
I_{(s)}=\int \sum_i\sum_k [(b_i-b_k)
(a_i-a_k)t+(a_i-a_k)^2 x]\ga i k\ga k i dx\eqno(5.9)
$$
\bigskip
\noindent
We emphasize the fact that the $t$--flow and the $s$--flow commute, so that
$$
\pois{\pd 1 I,\pd {(s)} I}-{\dd \pd {(s)} I\over \dd t}=0.
$$
By substituting:
$$
\int [Tr(\nabla_{_{u}}I_1[\nabla_{_{u}}I_{(s)},{d\over dx}+u])+{1\over
2}Tr(uv)]dx=0.  
$$
Then  there exists a function $Q_{_{(t)}}(x,t,u,v)$ such that
$$
Tr(-v[-tv-xu,{d\over dx}+u]+{uv\over 2})=-{d\over dx} Q_{_{(t)}}.
$$
By direct calculation (see Appendix 5.A) we obtain
$$
Q_{_{(t)}}={1\over 2}Tr(xuv+tv^2)
= {1\over 2}\sum_{i,j}
[(a_{_{j}}-a_{_{i}})(b_{_{i}}-b_{_{j}})x+(b_{_{j}}-b_{_{i}})^2t]\ga i j
\ga j i
$$

As in the previous examples, $Q_{_{(t)}}$ is the Hamiltonian for the reduced
$t$--flow. We  now describe this flow.

We start by rewriting the system
$$
\cases{u_t-v_x-[u,v]&=\ 0\cr
tu_t+xu_x+u&=\ 0\cr}\eqno{(5.10)}
$$
in terms of $\gamma$, i.e. we solve 
$$
[\l _t,a]\ =\ [\l _x,b]+[[\l,a],[\l,b]]
$$
under the condition
$$
\l _x=-{t\over x}\l _t-{1\over x}\l.
$$
This gives
$$
[\l _t,ax+tb]+[\l,b]=[[\l,ax],[\l,b]]
$$
but, because of the commutativity of $b$ with itself,
$$
{d\over dt}[\l,ax+tb]=[[\l,ax+bt],[\l,b]].\eqno{(5.11)}
$$
\bigskip

\noindent Then we identify $\S_s$ with the space of matrices 
$$q=[\l,
ax+bt]=ux+vt,
$$
so that $u=\ad_{_{(ax+bt)}}\ad_{_{a}}^{-1} q$, 
$v=\ad_{_{(ax+bt)}}\ad_{_{b}}^{-1} q$, or, in terms of the matrix elements: 
$$
q_{_{ij}}=[(a_{_{j}}-a_{_{i}})x+(b_{_{j}}-b_{_{i}})t]\ga i j
$$
On $\S$ the  equation (5.11) has the  Lax form
$$
q_t=[q,v]=[q,\ad_{_{(ax+bt)}}\ad_{_{b}}^{-1} q]
$$
with the Hamiltonian function
$$
H_{_{(t)}}= {1\over 2}\tr (q v)={1\over 2}\tr (x u v+t v^2).\eqno{(5.12)}
$$
This coincides with $Q_{_{(t)}}$.
\bigskip
\noindent
One may change the role of $x$ and $t$. This means that one considers
the system (5.1) on the space of functions $v(t)$:
$$
v_x=[\nabla_{_{v}}I_0,{d\over dt}+v],
$$
where $I_0=\int H_0 dt=-{1\over 2}\int \tr (uv)\  dt$ 
and one integrates in the
variable $t$. The scaling (5.8) can be read as an Hamiltonian equation
$$
{d v\over ds}=[\nabla_{_{v}}I_{(s)},{d\over dt}+v]=0,\eqno(5.13)
$$
where
$$
I_{(s)}=\int L dt=-{1\over 2}\int \tr\ (xuv+tv^2)\ dt\eqno(5.14)
$$
and
$$
\nabla_{_{v}} I_{(s)}=-tv-xu\eqno{(5.15)}
$$

Commutativity of the flows is equivalent to
$$
\pois{\pd 0 I,\pd {(s)} I}-{\dd \pd {(s)} I\over \dd x}=0;
$$
in our case
$$
\int [\tr(\nabla_{_{v}}I_0[\nabla_{_{v}}I_{(s)},{d\over dt}+v])+{1\over
2}Tr(uv)]dt=0.  $$
Then  there exists a function $Q_{_{(x)}}(x,t,u,v)$ such that
$$
\tr(-u[-tv-xu,{d\over dt}+v]+{uv\over 2})=-{d\over dt} Q_{_{(x)}}.
$$
By direct calculation (see Appendix 5.A) we obtain
$$
Q_{_{(x)}}={1\over 2}Tr(tuv+xu^2),\eqno{(5.16)}
$$
in components:
$$
Q_{_{(x)}}= {1\over 2}\sum_{i,j}
[(a_{_{j}}-a_{_{i}})(b_{_{i}}-b_{_{j}})t+(a_{_{j}}-a_{_{i}})^2x]\ga i j
\ga j i
$$
\bigskip
Now we study the $x$--flow on the reduced manifold defined by the scaling
equation :

the system (5.10) gives
$$
\cases{u_t-v_x-[u,v]&=\ 0\cr
tv_t+xv_x+v&=\ 0.\cr}\eqno{(5.17)}
$$
In terms of $\gamma$ this becomes
$$
{d\over dx}[\l,ax+tb]=[[\l,ax+bt],[\l,a]],
$$
that is a  Lax equation on $\S_s$:
$$
q_x=[q,u]=[q,\ad_{_{(ax+bt)}}\ad_{_{a}}^{-1} q]\eqno{(5.18)}
$$
with Hamiltonian function
$$
H_{_{(x)}}= {1\over 2}\tr (q u)={1\over 2}\tr (x u ^2+t  u v).\eqno{(5.19)}
$$
This coincides with $Q_{_{(x)}}$.

\bigskip

\noindent In fact one can rewrite the scaling as a
zero--curvature equation in two ways: 
$$
{d u\over ds}=q_x+[u,q]=0\eqno{(5.20)}
$$
and
$$
{d u\over ds}=q_t+[v,q]=0.\eqno{(5.21)}
$$

Therefore one may rewrite them  in terms of $q$ as
$$
q_x=[q,ad^{-1}_{_{(ax+bt)}} ad_{_{a}} q]
$$
and
$$
q_t=[q,ad^{-1}_{_{(ax+bt)}} ad_{_{b}} q].
$$
\bigskip
\noindent{\scap 5.2 Commuting time--dependent Hamiltonian flows on $\so(n)$}
\bigskip
\noindent We can do exactly the same using the coordinates 
$$
t_i=xa_i+tb_i,
$$
 and the
corresponding derivatives ${d\over dt_{_{i}}}$, with 
$${d\over dx}=\sum_i a_i 
{d\over dt_{_{i}}}
$$
and  
$$
{d\over dt}=\sum_i b_i {d\over dt_{_{i}}}.
$$

The starting equation  is now
$$
{d\over dt_{k}} u_i-{d\over dt_{_{i}}} u_k-[u_i,u_k]=0\eqno{(5.22)}
$$
where 
$$
u_i=[\l,E_i]\qquad (u_i)_{kl}=\ga k l\delta_{ik}-\ga k l\delta_{il}
$$
and $(E_i)_{kl}=\delta_{ik}\delta_{kl}$.
We impose the scaling
$$
{d\over ds} u_k=\sum_i t_i {d\over dt_{_{i}}} u_k\ +\ u_k=0\eqno{(5.23)}
$$

For every $k$ one can define, on the space 
$\F_k$ of functionals
$$
F=\int f(t,u,{du\over dt},....,{d^m u\over dt^m_k})\ dt
$$
with
$$
{d\over d\epsilon}F(u_k+\epsilon w)\vert_{_{\epsilon=0}}=(\nabla_{_{u_k}}F,w),
$$
a Poisson structure $P^{^{(k)}}$ with the Poisson
bracket
$$	
\{F,G\}(u_k)=(\nabla_{_{u_k}}F,[\nabla_{_{u_k}}G,{d\over dt_{k}}
+u_k])
$$
The $n$-waves equation (5.22) is Hamiltonian w.r.t.  the Poisson structure
 $P^{^{(k)}}$ in $\F_k$: 

$$
{d\over dt_{_{i}}} u_k=[\nabla_{_{u_k}}I_{i},{d\over
dt_{k}}+u_k]=[-u_i,{d\over dt_{k}}+u_k],\eqno(5.24) $$
where
$$
I_{i}=\int L_{i} dt=-{1\over 2}\int \tr\ (u_i u_k)\ dt_k.\eqno{(5.25)}
$$

On $\F_k$ we  can reduce to P VI equation imposing the scaling (5.23),
which admits the Hamiltonian form
$$
{d\over ds} u_k=[\nabla_{_{u_k}}I_{(s)},{d\over dt_{k}}+u_k]=
[-\sum_j t_j u_j,{d\over dt_{k}}+u_k]=0\eqno{(5.26)}
$$
where
$$
I_{(s)}=\int L dt_k=-{1\over 2}\int \tr\sum_j t_j u_j u_k dt_k
=-{1\over 2}\int \tr(\sum_{j\not= k} t_j u_j u_k +t_k u_k^2) dt_k\eqno{(5.27)}
$$

The commutativity of the flows  is equivalent to
$$
\pois{I_{i},I_{(s)}}-{\dd I_{(s)}\over \dd t_i}=0;
$$
in our case
$$
\int [\tr(\nabla_{_{u_k}}I_{i}[\nabla_{_{u_k}}I_{(s)},
{d\over dt_{k
}}+u_k])+{1\over
2}Tr(u_i u_k)]dt_k=-\int \dd_k Q_{_{(i)}} dt_k.  
$$
By direct calculation (following the scheme in Appendix 5.A) we obtain
$$
Q_{_{(i)}}={1\over 2}\tr \sum_jt_j u_j u_i\ =\ \sum_j(t_i-t_j)\ga i j\ga j
i\eqno{(5.28)} $$

\bigskip
The scaling equation defines the submanifold $\S_s$. One can consider on $\S_s$
the system of coordinates given by the matrix elements of $q$:
$$
q=[\l,\sum_j t_j E_j]=[\l,U],
$$
where $U$ is the diagonal matrix diag$(t_1,\dots,t_n)$; explicitly
$$
q_{ij}=(t_j-t_i)\ga i j.\eqno {(5.29)}
$$
As in the previous cases, $Q_{_{(i)}}$ is the Hamiltonian for the $t_i$-flow on
the reduced manifold.

Starting now from $\F_i$, $i\not= k$,
we can reduce on the same submanifold $\S$ and construct the Hamiltonian
function  $Q_{_{(k)}}$. 

Indeed, the scaling  (5.23) for every $k$ produces on $\S_s$ the Lax
equation 
$$
q_k=[q,u_k],\eqno{(5.30)}
$$
with Hamiltonian functions 
$$
H_k= {1\over 2}\tr (q u_k)={1\over 2}\tr \sum_jt_j u_j u_k=
{1\over 2}\sum_{j\not= k}{q_{jk}q_{kj}\over t_k-t_j}.\eqno{(5.31)}
$$
These coincide with the $Q_{_{(k)}}$ constructed above.
Observing that $\l=\ad^{-1}_{_{U}}q$ one can rewrite  
$$
u_k=\ad_{_{E_k}}\ad^{-1}_{_{U}}q.
$$
In the case $q^{^{T}}=-q$ eqs. (5.30) are the Monodromy Preserving 
Deformation equations for the linear differential operator 
$$
\Lambda={d\over d\lambda}-U-{q\over\lambda}
$$
that give Painlev\'e VI, for $n=3$, and the higher--order analogues, 
for $n>3$.
\bigskip 
\noindent {\bf Remark}: The first integrals of the MPDE (5.30) are given
by the monodromy data of the operator $\Lambda$. The Poisson bracket on the space of the monodromy data has been computed in [Ug].
\bigskip
\bigskip
\bigskip

\noindent {\scap 5.A Appendix}
\bigskip
Here we present the explicit calculations giving rise to equation (5.16):
let us consider the following explicit expressions:
$$
\eqalignno{I_t&=-{1\over 2}\int \tr\ (uv)dx\cr
\nabla I_t&=-v\cr
I_{(s)}&=-{1\over 2}\int \tr\ (xu^2+tuv)dx\cr
\nabla I_{(s)}&=-tv-xu\cr
\pois{I_t,I_{(s)}}&=(-v,tu_t+xu_x+u)=\cr
&=-\int \tr\ 
(tvv_x+tv[u,v]+xvu_x+uv)dx\cr
{\dd \pd {(s)} I\over \dd t}&=I_t=-{1\over 2}\int \tr\ (uv).\cr 
\pois{I_t,I_{(s)}}-{\dd \pd {(s)} I\over \dd t}&=
-\int \tr\ 
(tvv_x+tv[u,v]+xvu_x+{1\over2}uv)dx&(a.1)\cr}
$$
In (a.1) the relations 
$$
\eqalign{\tr \ v\ [u,v]&=0\cr
\tr\ (vv_x)&={1\over
2}{d\over dx} \tr\ (v^2)\cr
\tr\ (xvu_x+{1\over2}uv)&={1\over2}{d\over dx}\tr\ (xuv)\cr}
$$ 
hold. In fact, in terms of $\ga i j$ one can write 
$$
\eqalign{\tr\ (xvu_x)
&=\sum_i\sum_k x(b_k-b_i)(a_i-a_k)\ga i k(\gamma_x)_{ki} =\cr
&=\sum_i\sum_k x(b_i-b_k)(a_k-a_i)(\gamma_x)_{ki}\ga i k=\cr
&=\tr\ (xv_xu),\cr}
$$

which implies
$$
{d\over dx}\tr\ (xuv)=2\tr\ (xvu_x)+\tr\ (uv).
$$
Then:
$$
\pois{I_t,I_{(s)}}-{\dd \pd {(s)} I\over \dd t}
=-{1\over 2}\int {d\over dx}\tr\ (x u v+tv^2)\ dx.
$$
\bigskip

\noindent{\scap Acknowledgements}
\bigskip
I wish to thank prof. Boris Dubrovin for his  guidance into the subject, for 
his precious suggestions and corrections and for the final reading of the manuscript.
My thanks also to Gregorio Falqui for many  useful discussions.
\bigskip
\noindent{\scap References}
\bigskip
\bigskip
\noindent [Al] S.I.Al'ber: {\it Investigation of equations of Korteweg--de
Vries type by the method of recurrence relations}, J.London.Math.Soc.(2) {\bf
19}, 467-480 (1979).

 \bigskip
\noindent [AS] M.J. Ablowitz, H. Segur: {\it Solitons and the Inverse 
Scattering Transform}, SIAM, Philadelphia (1981).
\bigskip

\noindent [BN] O.I. Bogoyavlenskii, S.P. Novikov: {\it The relationship
between Hamiltonian formalism of stationary and nonstationary problems},
Funct. Anal. Appl. {\bf 10} n.1 (1976), 9--13.
\bigskip
\noindent [CD] F. Calogero, A. Degasperis: {\it Spectral transform and 
solitons},  North--Holland Publishing Company (1982).
\bigskip

\noindent [D1] B. Dubrovin:
{\it Painlev\'e transcendents in two--dimensional topological field theory},
Proceedings of 1996 Carg\`ese
summer school `` Painlev\'e Transcendents: One Century Later".





\bigskip
\noindent [Mo] O.I.Mokhov:
{\it On the Hamiltonian property of an arbitrary evolution system on the set
of stationary points of its integral},  Communication of the Moscow Math. Soc.
{\bf }, 133-134,  (1983). 
\medskip
\quad O.I.Mokhov:
{\it The Hamiltonian property of an  evolutionary flow on the set
of stationary points of its integral},  Math.Ussr Izvestiya {\bf 31}, 657-664, 
(1988).
  
\bigskip \noindent [MW] L.J.Mason, N.M.J. Woodhouse: {\it
Integrability,  Self--Duality and Twistor Theory},  L.M.S. Monographs New Series
{\bf 15}
 Oxford Univ. Press, (1996).
\bigskip
\noindent [O] K. Okamoto: {\it The Painlev\'e equations and the 
Dynkin diagrams}, in ''Painlev\'e Transcendents", Eds. D. Levi
and P. Winternitz, Plenum  Press, New York (1984).
\bigskip
\noindent [SM] A.B. Shabat, A.V. Mikhailov:
{\it  Symmetries -- Test of integrability}, in ''Important Developments 
in Soliton Theory", Eds. A.S. Fokas, V.E. Zakharov, Springer Series in 
Nonlinear Dynamics.

\bigskip
\noindent [Ug] M.Ugaglia:
{\it On  a Poisson structure on the space of Stokes matrices},  IMNR, 9,
473-493 (1999).
\bigskip 

\bye